\documentclass[letterpaper]{article}%
\usepackage{graphicx}
\usepackage{amssymb}
\usepackage{amsmath}
\usepackage{amsfonts}
\usepackage{sw20bams}
\usepackage{comment}
\usepackage{hyperref}
\usepackage{url}
\usepackage{aurl}
\providecommand{\U}[1]{\protect\rule{.1in}{.1in}}

\ifx\pdfoutput\relax\let\pdfoutput=\undefined\fi
\newcount\msipdfoutput
\ifx\pdfoutput\undefined\else
\ifcase\pdfoutput\else
\ifx\paperwidth\undefined\else
\ifdim\paperheight=0pt\relax\else\pdfpageheight\paperheight\fi
\ifdim\paperwidth=0pt\relax\else\pdfpagewidth\paperwidth\fi
\fi\fi\fi
\pdfoutput=1
\begin{document}

\title{Fast Ramanujan--type Series for Logarithms. Part II.}
\author{Jorge Zuniga\thanks{$\ \ $Independent Researcher \texttt{jorge.zuniga.hansen.at.gmail.com}}}
\date{June 2025\footnote{$\ \ $Keywords: NT number theory, Sequences and Series,
Algorithms, Binary Splitting, LLL method, Mixed Integer Linear Programming,
Monte Carlo method}\vspace*{-4pt}}
\maketitle

\begin{abstract}
This work extends the results of the preprint \textit{Fast Ramanujan-type
Series for Logarithms, Part I}, $\small{\text{arXiv:}2506.08245}$, which introduced single
hypergeometric type identities for the efficient computing of $\log(p)$, where
$p\in%
\mathbb{Z}
_{>1}$. Novel formulas for arctangents and very fast
evaluation methods of multiple logarithms are introduced. Building upon a $\mathcal{O}((p-1)^{6})$
series asymptotic approximation for $\log(p)$ as $p\rightarrow1$, Eq.(45)
\textit{Part I}, 
formulas for computing $n$ simultaneous logarithms are developed.
These formulas are derived by solving an integer programming problem to
identify optimal variable values within a finite lattice $\subset%
\mathbb{Z}
^{n}.$ This approach yields linear combinations of hypergeometric series that provide:
{\normalsize (i)} the most efficient formulas for single logarithms of natural
numbers (some of them were tested to compute more than $10^{11}$ decimal digits) and {\normalsize (ii)} the fastest known formulas for
multi-valued logarithms of $n$ selected integers in $%
\mathbb{Z}
_{>1}$.

\end{abstract}

\section{Introduction\smallskip}

Efficient algorithms for the simultaneous computation of multiple logarithms
(multi-valuation) are crucial in various numerical analysis applications. For
example, van der Hoeven and Johansson's meta-algorithm (1) in \cite{FRED}
utilizes multi-valued logarithms for high-precision computation of the
exponential function $e^{x}$, a technique applied to other elementary
functions \cite{FREDJ}\cite{SCH01}\cite{SCH02} and also to Riemann's
$\zeta(z)$ function and its derivatives \cite{JOHAN}.

This article extends these applications by presenting novel formulas for
arctangents and multi-valued logarithms, leading to highly efficient
algorithms for multiprecision calculation of single logarithms. The
approach relies on Eq.(45) from Part\daurl{:I}{https://arxiv.org/pdf/2506.08245}\aurl{:I}{}, yielding rapid
hypergeometric series for $\log(x)$ and a$\tanh(x)$ suitable for binary
splitting computations.

In this work formulas for $n$ simultaneous logarithms --see \cite{FREDJ}, Algorithm 2, Step (1)-- are built by linearly
independent combinations of $n$ fast hypergeometric series. The parameters of
these series are determined by solving an integer non-linear programming
problem whose set of feasible solutions, with optimal variable values within a
lattice $\subset%
\mathbb{Z}
^{n}$, is sorted by computing cost. These solutions are searched within and on
the edges of polyhedra defined by a bounded linear objective function and
bilaterally bounded integer variables.

PARI-GP \cite{PARIGP} scripts (See Section 7, PARI-GP Modules) were coded to obtain this set of solutions
using three distinct methods: (i) exhaustive brute-force iteration, (ii) external mixed-integer programming (MIP) solvers (SCIP\texttrademark \thinspace\cite{SCIP}/COPT\texttrademark \thinspace\cite{COPT}), 
and (iii) a randomized Monte Carlo algorithm that
employs the LLL lattice reduction algorithm (via PARI's \textit{lindep}
function) with fine-tuned internal precision to capture feasible solutions,
each associated to its own hypergeometric series. Methods (i) and (ii) are suitable for small $n<8$ and (iii) for medium  $n<50$. For all methods, the
top-ranked $n$ fastest series are selected based on cost, provided their summand's
polynomial integer coefficients do not exceed a given bit-size limit. The final formulas yield arbitrary precision constants and several of them were tested to compute $3.4\cdot10^{11}$ decimal places for some single logarithms.

\section{Logarithms and arctangents fast hypergeometric series\medskip}

Ramanujan type hypergeometric series of the type ${}_{4}F_{3}$ for $\log(x), $
atanh$(x)$ and a$\tan(x)$ with convergence in $x\in\mathcal{D}\subseteq\mathbb{C}$, but defined only for $x\in\mathcal{D}\cap\mathbb{Q}$,    
are introduced as follows,\vspace*{-1pt}

\subsection{log(\textit{x}) series}

For $u,v\in%
\mathbb{Z}
_{>0}$ such that $u\neq v$ and $\gcd(u,v)=1,$ Eq.(45) Part\daurl{:I}{https://arxiv.org/pdf/2506.08245}\aurl{:I}{} with
$p=x=\tfrac{u}{v}$ and convergence in $\mathcal{D}_1=\{x|x\in%
\mathbb{Q}
\wedge|x-7|<4\sqrt{3}\}$ can be written as\vspace*{-1pt}%
\begin{equation}
\log\left(  \dfrac{u}{v}\right)  =-\frac{1}{12}\cdot\frac{u-v}{u^{2}%
v^{2}(u+v)}\cdot\sum_{n=0}^{\infty}\frac{P_{1}(n,u,v)}{(6n+1)(6n+5)}%
\cdot\left(  \frac{(u-v)^{6}}{108\,u^{2}v^{2}(u+v)^{2}}\right)  ^{\!\!n}\cdot%
\begin{bmatrix}%
\genfrac{}{}{0pt}{0}{{}}{{}}%
1 & \frac{1}{2}%
\genfrac{}{}{0pt}{0}{{}}{{}}%
\\%
\genfrac{}{}{0pt}{0}{{}}{{}}%
\frac{1}{6} & \frac{5}{6}%
\genfrac{}{}{0pt}{0}{{}}{{}}%
\end{bmatrix}
_{n}\tag{1}\label{1}%
\end{equation}\vspace{-2pt}
\begin{equation}
\text{ \ \ \ \ \ \ \ \ \ }=\sum_{n=0}^{\infty}\frac{-(u-v)\cdot P_{1}%
(n,u,v)}{12\,u^{2}v^{2}(u+v)(6n+1)(6n+5)}\cdot\prod_{k=0}^{n-1}\frac
{(u-v)^{6}(k+1)(2k+1)}{6\,u^{2}v^{2}(u+v)^{2}(6k+1)(6k+5)}\tag{2}\label{2}%
\end{equation}
\[
P_{1}(n,u,v)=2\,(u^{2}-14uv+v^{2})(u^{2}+4uv+v^{2})\cdot n+u^{4}%
-14u^{3}v-94u^{2}v^{2}-14uv^{3}+v^{4}\vspace*{4pt}%
\]
where the $\Sigma\,\Pi$ form matches FLINT's \cite{FLINT0} binary splitting
hypergeometric convention.

\subsection{atanh(\textit{x}) series}

By setting $u=t+d,$ $v=t-d$ a series for the atanh function is obtained, Eq.(\ref{8}), with
$x=\frac{d}{t}$ and region of convergence $\mathcal{D}_2=\{x|x\in%
\mathbb{Q}
\wedge|\,x\,|<\sqrt{3}/2\}$
\begin{equation}
\text{atanh}\left(  \frac{d}{t}\right)  =\frac{1}{3}\cdot\frac{d}%
{t\,(t^{2}-d^{2})^{2}}\cdot\sum_{n=0}^{\infty}\frac{P_{2}(n,d,t)}%
{(6n+1)(6n+5)}\cdot\left(  \frac{4\,d^{6}}{27\,t\,^{2}(t^{2}-d^{2})^{2}%
}\right)  ^{\!\!n}\cdot%
\begin{bmatrix}%
\genfrac{}{}{0pt}{0}{{}}{{}}%
1 & \frac{1}{2}%
\genfrac{}{}{0pt}{0}{{}}{{}}%
\\%
\genfrac{}{}{0pt}{0}{{}}{{}}%
\frac{1}{6} & \frac{5}{6}%
\genfrac{}{}{0pt}{0}{{}}{{}}%
\end{bmatrix}
_{n}\tag{3}\label{3a}%
\end{equation}\vspace{-2pt}
\begin{equation}
\text{ \ \ \ \ \ \ \ \ \ \ \ }=\sum_{n=0}^{\infty}\frac{d\cdot P_{2}%
(n,d,t)}{3\,t\,(t^{2}-d^{2})^{2}(6n+1)(6n+5)}\cdot\prod_{k=0}^{n-1}%
\frac{8\,d^{6}(k+1)(2k+1)}{3t^{2}(t^{2}-d^{2})^{2}(6k+1)(6k+5)}\tag{4}%
\label{4a}%
\end{equation}
\[
P_{2}(n,d,t)=2\,\,(3t^{2}-4d^{2})(3t^{2}-d^{2})\cdot n+15t^{4}-25t^{2}%
d^{2}+8d^{4}%
\]

\subsection{atan(\textit{x}) series}

replacing $d$ by $i\cdot d,$ an hypergeometric alternating series for the atan
function is achieved with $x=\frac{d}{t}$ converging in\ $\mathcal{D}_3=\{x|x\in%
\mathbb{Q}
\wedge|\,x\,|<2.90578...\}$\vspace*{-1pt}
\begin{equation}
\text{atan}\left(  \frac{d}{t}\right)  =\frac{1}{3}\cdot\frac{d}%
{t\,(t^{2}+d^{2})^{2}}\cdot\sum_{n=0}^{\infty}\frac{P_{3}(n,d,t)}%
{(6n+1)(6n+5)}\cdot\left(  \frac{-4\,d^{6}}{27\,t\,^{2}(t^{2}+d^{2})^{2}%
}\right)  ^{\!\!n}\cdot%
\begin{bmatrix}%
\genfrac{}{}{0pt}{0}{{}}{{}}%
1 & \frac{1}{2}%
\genfrac{}{}{0pt}{0}{{}}{{}}%
\\%
\genfrac{}{}{0pt}{0}{{}}{{}}%
\frac{1}{6} & \frac{5}{6}%
\genfrac{}{}{0pt}{0}{{}}{{}}%
\end{bmatrix}
_{n}\tag{5}\label{5}%
\end{equation}%
\begin{equation}
\text{ \ \ \ \ \ \ \ \ \ }=\sum_{n=0}^{\infty}\frac{d\cdot P_{3}%
(n,d,t)}{3\,t\,(t^{2}+d^{2})^{2}(6n+1)(6n+5)}\cdot\prod_{k=0}^{n-1}%
\frac{-8\,d^{6}(k+1)(2k+1)}{3t^{2}(t^{2}+d^{2})^{2}(6k+1)(6k+5)}%
\tag{6}\label{6}%
\end{equation}%
\[
P_{3}(n,d,t)=2\,\,(3t^{2}+4d^{2})(3t^{2}+d^{2})\cdot n+15t^{4}+25t^{2}%
d^{2}+8d^{4}\vspace*{3pt}%
\]

Depending on $u,v,d,t$ these series can be very fast formulas specially fitted
to apply the binary splitting method for multiprecision computing. For large
$n,$ $\log(x)$ series terms are $\mathcal{O}((x-1)^{6n})$ as $x\rightarrow1$
while for atanh$(x)$ and atan$(x)$\ they are $\mathcal{O}(x^{6n})$ as
$x\rightarrow0.$

\subsubsection{Example 1}

The recurrence\vspace*{-1pt}%
\begin{equation}
\log(n+1)=\log(n)+2\,\text{atanh}\left(  \tfrac{1}{2n+1}\right)  ,\text{
\ }n\in%
\mathbb{N}%
\tag{7}\label{7}%
\end{equation}
can be computed using series Eqs.(\ref{3a}--\ref{4a}) if the previous
logarithm is known. For $n=14$ the binary splitting cost --Eq.(7) Part\daurl{:I}{https://arxiv.org/pdf/2506.08245}\aurl{:I}{}-- 
is about the same (and so it is as fast) as Chudnovky's $\pi$ formula
\cite{CHUD}. Greater values of $n$ produce much faster algorithms. This should
be a standard method to compute a sequence of logarithms of consecutive
natural numbers. It requires to compute by binary splitting only one fast
series per logarithm.

\subsubsection{Remark 1}

Note that Eqs.(\ref{1}--\ref{2}) and Eqs.(\ref{3a}--\ref{4a}) are two
equivalent series. $\log(x)$ and a$\tanh(x)$ are linked by these elementary
identities\vspace*{0pt}
\begin{equation}%
\begin{tabular}
[c]{lcl}%
$\log\left(  \dfrac{u}{v}\right)  $ & $=$ & 2\thinspace$\text{atanh}\left(
\dfrac{u-v}{u+v}\right)  $\\
&  & \\
$\text{atanh}\left(  \dfrac{d}{t}\right)  $ & $=$ & $\dfrac{1}{2}\log\left(
\dfrac{t+d}{t-d}\right)  $%
\end{tabular}
\tag{8}\label{8}%
\end{equation}
where $u,v,d,t$ ratios in argument $x$ must give $x$ belonging to the series
convergence regions indicated above, namely $\mathcal{D}_1\cap\{x|$ $x>0\}$ for
$\log(x)$ and $\mathcal{D}_2\cap\{x|$ $|x|<1\}$ for atanh$\left(  x\right)  $.

\subsubsection{Remark 2}

Eqs.(\ref{1}--\ref{2}), Eqs.(\ref{3a}--\ref{4a}) and Eqs.(\ref{5}--\ref{6})
provide high performance identities for multiple precision computing of
$\log(x)$, a$\tanh(x)$ and a$\tan(x)$ functions with $x\in%
\mathbb{Q}
$. These series have several applications including the fast computing of
$\exp(z),$ $\log(z)$ and elementary trigonometric functions for $z\in
\mathcal{D}_{z}\subseteq%
\mathbb{C}
$ in some region $\mathcal{D}_{z}\ $\cite{FREDJ}.

\subsubsection{Example 2}

The hypergeometric series origin of these formulas allows us to measure the
efficiency through the binary splitting cost, Eq.(7) Part\daurl{:I}{https://arxiv.org/pdf/2506.08245}\aurl{:I}{}, enabling to compare different algorithms. For example \cite{FREDJ} --\textit{Table\ }1--\ that is based on Luca-Najman table \cite{CMT}\cite{FL1}\cite{FL2}, has\ Machin--like formulas for multivaluation of
logarithms of the first $n$ primes.\ For\ $n=2,3,4,5,6,$\ \textit{Table\ }%
1\ data\ applied\ on\ Eqs.(\ref{3a}{\small --}\ref{4a}) give integer
coefficients of bit-size 32, 64, 96, 96 and 128 and total binary
splitting costs 1.014, 0.919, 0.729, 0.823 and 0.837 respectively. However the
optimization method introduced here (see next section) finds solutions
 giving multi-valuated logarithm hypergeometric formulas with binary splitting costs
1.014, 0.819, 0.694, 0.779 and 0.756 for the same primes and same bit--size of input
integers respectively. These are pretty better than Machin formulas resulting in faster
algorithms. In addition, Machin can be outperformed much further if input bit
sizes \textit{b} are relaxed and larger \textit{b} values are allowed.

\subsection{y-cruncher binary splitting series}

Eq.(\ref{1}), Eq.(\ref{3a}), Eq.(\ref{5})\ can be also written using
y-cruncher's \cite{YCRUNCHER} binary splitting convention for $x\in%
\mathbb{Q}
$ as

\subsubsection{log(\textit{x}) series}%

\begin{equation}
\log\left(  \frac{u}{v}\right)  =-\frac{u^{2}-v^{2}}{2}\cdot\sum_{n=1}%
^{\infty}\frac{\mathcal{P}_{1}(n,u,v)}{(u-v)^{6}\,n\,(2n-1)}\cdot\prod
_{k=1}^{n}\frac{(u-v)^{6}\,k\,(2k-1)}{6\,u^{2}v^{2}(u+v)^{2}\,(6k-1)(6k-5)}%
\tag{9}\label{9}%
\end{equation}%
\[
\mathcal{P}_{1}(n,u,v)=2\,(u^{2}-14uv+v^{2})(u^{2}+4uv+v^{2})\cdot
n-(u+v)^{2}(u^{2}-8uv+v^{2})\vspace*{4pt}%
\]

\subsubsection{atanh(\textit{x}) series}%

\begin{equation}
\text{atanh}\left(  \frac{d}{t}\right)  =d\,t\cdot\sum_{n=1}^{\infty}%
\frac{\mathcal{P}_{2}(n,d,t)}{8\,d^{\,6}\,n(2n-1)}\cdot\prod_{k=1}^{n}%
\frac{8\,d^{\,6}\,k(2k-1)}{3\,t^{2}(t^{2}-d^{2})^{2}(6k-1)(6k-5)}%
\tag{10}\label{10}%
\end{equation}%
\[
\mathcal{P}_{2}(n,d,t)=2\,\,(3t^{2}-4d^{2})(3t^{2}-d^{2})\cdot n-t^{2}%
(3t^{2}-5d^{2})\vspace*{4pt}%
\]

\subsubsection{atan(\textit{x}) series}%

\begin{equation}
\text{atan}\left(  \frac{d}{t}\right)  =d\,t\cdot\sum_{n=1}^{\infty}%
\frac{\mathcal{P}_{3}(n,d,t)}{-8\,d^{\,6}\,n(2n-1)}\cdot\prod_{k=1}^{n}%
\frac{-8\,d^{\,6}\,k(2k-1)}{3\,t^{2}(t^{2}+d^{2})^{2}(6k-1)(6k-5)}%
\tag{11}\label{11}%
\end{equation}%
\[
\mathcal{P}_{3}(n,d,t)=2\,\,(3t^{2}+4d^{2})(3t^{2}+d^{2})\cdot n-t^{2}%
(3t^{2}+5d^{2})
\]

\section{Algorithms for multi-valuation of logarithms\medskip}

This section introduces new algorithms for the simultaneous evaluation of
multiple logarithms, derived from solving an integer non-linear optimization
problem, and based on Eq.(45) Part\daurl{:I}{https://arxiv.org/pdf/2506.08245}\aurl{:I}{}, Eqs.(\ref{1}--\ref{2}) and
Eq.(\ref{9}).

\subsection{An integer non-linear programming problem}

Let us consider the vector $\boldsymbol{q}=$ $[\,p_{1},p_{2},...p_{n}\,]$ of
$n$ different integers in $%
\mathbb{Z}
_{>1}$ with exactly $n$ common prime factors (i.e. logarithms are linearly
independent). Generality is not lost if $p_{j},\,j=1,...n$ are indeed $n$
different primes and, unless it is stated otherwise, this is so in what
follows. Consider also a vector $\boldsymbol{x}=[$\thinspace$x_{1}%
,x_{2},...x_{n}\,]$ $\in%
\mathbb{Z}
^{n}$ of integer power exponents not all equal to 0 (to exclude the trivial case)
such that for $p\in\mathbb{Q}$ and $p>0$,\vspace*{-3pt}%
\begin{equation}
p=\prod_{j=1}^{n}p_{j}^{x_{j}}\text{\ \ }\Longrightarrow\text{\ \ }\log
p=\sum_{j=1}^{n}x_{j}\cdot\log p_{j}=\sum_{j=1}^{n}x_{j}\cdot c_{j}%
=\boldsymbol{x}\cdot\boldsymbol{c}=\mathcal{S}\left(  p\right)  \vspace
*{-2pt}\tag{12}\label{12}%
\end{equation}
holds, where $\boldsymbol{c}=[\,c_{1},c_{2},...c_{n}\,]^{\text{{\tiny T}}}=[\,\log
p_{1},\log p_{2},...\log p_{n}\,]^{\text{{\tiny T}}}$ is the column vector of
logarithms and $\mathcal{S}(p)$ is the RHS series in Eq.(45) Part\daurl{:I}{https://arxiv.org/pdf/2506.08245}\aurl{:I}{},
indicating that to get a fast $\log\hspace*{1pt}p$ series, $p$ should be as
close to 1 as possible. So $\varepsilon=|\log p\,|$ for some small positive
$\varepsilon,$ but the absolute value can be ignored because of the exponents
symmetry between $p$ and $p^{-1}$. Thus the searching of the fastest series
for multi-valuation of logarithms can be settled, at first glance, as an
integer linear programming (ILP) problem where $\varepsilon>0$ is a linear
objective function to be minimized\vspace*{-11pt},
\begin{equation}%
\begin{array}
[c]{c}%
\genfrac{}{}{0pt}{0}{\genfrac{}{}{0pt}{1}{{}}{{}}}{{}}%
\min\limits_{\text{ }\boldsymbol{x\in%
\mathbb{Z}
}^{n}}\,\varepsilon\\
\text{{\small \textit{s.t.}}}%
\genfrac{}{}{0pt}{0}{\genfrac{}{}{0pt}{1}{{}}{{}}}{{}}%
\varepsilon=\boldsymbol{x}\cdot\boldsymbol{c}%
\end{array}
\vspace*{-4pt}\tag{13}\label{13}%
\end{equation}
where $\boldsymbol{x}=[\,x_{1},x_{2},...x_{n}\,]\neq\lbrack\,0,0,...0\,]$ is
the vector of integer decision variables and $\boldsymbol{c}$ is the linear
programming column vector of costs. This problem, thus posed, has infinite
solutions, so restrictions to $\boldsymbol{x}$ must be imposed to confine the
solutions search space. This is done limiting the sizes of numerator and
denominator of $p$ by setting bounds $X_{j}$ on each $x_{j},$ $\,j=1,...n$
as\vspace*{-7pt}
\begin{equation}
-X_{j}\leq x_{j}\leq X_{j}\text{ \ with \ \ }X_{j}=\left\lceil \tfrac{M}%
{\log_{2}(p_{j})}\right\rceil \tag{14}\label{14}%
\end{equation}
for some fixed positive integer $M$ such that $2^{-M}\leq p_{j}^{x_{j}}%
\leq2^{M}.$ These bilateral bounds on $\boldsymbol{x}$ confine the solution
space to $n$ dimensional closed polyhedra.\ Optionally the objective function
can also be constrained by above if there is enough information to do it. So%
\begin{equation}
\boldsymbol{x}\cdot\boldsymbol{c}\leq\mathcal{E}\tag{15}\label{15}%
\end{equation}
for some given $\mathcal{E}$. As it has been raised Eqs.(\ref{13}%
--\ref{14}--\ref{15}) belong to an ILP problem class and can be solved relying
on branch-and-bound and branch-and-cut\ algorithms \cite{CORNELL} that have
been implemented in some known MILP softwares (mixed integer linear
programming where some of the decision variables can be continuous)
\cite{MILP}.

\subsubsection{A non-linear constraint}

By enforcing the numerator/denominator polynomial's coefficients $\mathcal{C}\in\mathbb{Z}$ in
series $\mathcal{S}(p)$ --see RHS of Eq.(\ref{2}) or Eq.(\ref{9})--
to be limited by some input bit size $b,$ this constraint\vspace*{-4pt}%
\begin{equation}
\left\lceil \log_{2}(\,\max|\,\mathcal{C}\,|\,)\right\rceil <b\vspace
*{2pt}\tag{16}\label{16}%
\end{equation}
converts the problem into a highly non-linear one. A strategy to address it is
to solve the ILP problem discarding the solution if Eq.(\ref{16}) is not
satisfied. In this case the next solution with greater $\varepsilon$ must be
searched by repeating this ILP solve and (brute force) verify process. A
feasible solution is found and accepted if this non-linear constraint holds.

\subsection{General setup, feasible solutions and multi-valuation identities}

The ILP problem Eqs.(\ref{13}--\ref{14}), with optionally Eq.(\ref{15}) and/or
Eq.(\ref{20}) further on, defines a finite sequence of feasible
solutions $\boldsymbol{\varepsilon}=[\,\varepsilon_{1},\varepsilon
_{2},...,\varepsilon_{m}\,]$ (also known as the solutions pool) with $m>n$
and
\begin{equation}
0<\varepsilon_{1}<\varepsilon_{2}<\varepsilon_{3}<...<\varepsilon_{m}%
\tag{17}\label{17}%
\end{equation}
where $\varepsilon_{1}$ is the global minimum, $\varepsilon_{2}$ is the second
best, $\varepsilon_{3}$ the third one and so on, so that
$\boldsymbol{\varepsilon}$ is written as a sorted vector by definition. This
sorting is equivalent to rank the resulting series by convergence rate, by
binary splitting cost or by computing speed (from faster to slower). For each
$\varepsilon_{i},\,i=1,2,...,m$ the vector $\boldsymbol{x}_{i}=[\,x_{i1}%
,x_{i2},...x_{in}\,]$ is the feasible solution corresponding to the $i^{th}$
row in a sorted $m\,${\scriptsize x}$\,n$\ matrix containing the ILP solutions
pool with the primes' exponents $x_{ij}$. For a logarithm vector
$\boldsymbol{c}$ of size $n$, this matrix is scanned selecting the first $n$
top ranked rows that give a full rank non-singular $n\,${\scriptsize x}%
$\,n$\ matrix and satisfy all constraints including non-linear ones
--Eq.(\ref{16})-- with $p=p^{(i)},$ $i=1,2,...,n$ and
\begin{equation}
p^{(i)}=\prod_{j=1}^{n}p_{j}^{x_{ij}}\text{\ \ }\Longrightarrow\text{\ \ }\log
p^{(i)}=\sum_{j=1}^{n}x_{ij}\cdot\log p_{j}=\sum_{j=1}^{n}x_{ij}\cdot
c_{j}=\boldsymbol{x}_{i}\cdot\boldsymbol{c}=\mathcal{S(}p^{(i)})\tag{18}%
\label{18}%
\end{equation}
where $\boldsymbol{x}_{i}$ are now the ordered selected rows belonging to this
non-singular solution matrix of integer values. Name this $n\,${\scriptsize x}%
$\,n$ matrix as ${}_{b}\boldsymbol{X}_{\boldsymbol{q}}$ to explicitly indicate
the dependence on $\boldsymbol{q}$, the vector of primes, and $b$ the target
formula's integer bit size limit shown in Eq.(\ref{16}). From this solution
matrix ${}_{b}\boldsymbol{X}_{\boldsymbol{q}}$ the column vector ${}%
_{b}\boldsymbol{S}_{\boldsymbol{q}}=\mathtt{[}\mathcal{\,S(}p^{(1)}),$
$\mathcal{S(}p^{(2)}),...\mathcal{S(}p^{(n)})\,\mathtt{]}^{\text{{\tiny T}}}%
$\thinspace of series is built taking $p^{(i)}$ from Eq.(\ref{18}) and
$\mathcal{S(}p)$ as the RHS series from Eq.(45) Part\daurl{:I}{https://arxiv.org/pdf/2506.08245}\aurl{:I}{}. 
Alternatively $_{b}\boldsymbol{S}_{\boldsymbol{q}}$ series elements are computed from
numerator $u_{i}$ and denominator $v_{i}$ in $p^{(i)}=u_{i}/v_{i}%
,\ i=1,2,...n$ separating $u_{i},v_{i}$ factors $p_{j}^{x_{ij}}$ according to
the signs of $x_{ij}$ in Eq.(\ref{18}) (See $\it{Appendix}$). $_{b}\boldsymbol{S}_{\boldsymbol{q}}$
series are then lifted up by applying Eqs.(\ref{1}--\ref{2}) or Eq.(\ref{9}).

\bigskip

\hspace*{-15pt}The final multi-valuation identity is obtained solving this
linear system of equations
\begin{equation}
{}{}_{b}\boldsymbol{X}_{\boldsymbol{q}}\cdot\boldsymbol{c}={}_{b}%
\boldsymbol{S}_{\boldsymbol{q}}\Longrightarrow\boldsymbol{c}={}_{b}%
\boldsymbol{X}_{\boldsymbol{q}}^{-1}\cdot{}_{b}\boldsymbol{S}_{\boldsymbol{q}%
}\tag{19}\label{19}%
\end{equation}
and each logarithm in $\boldsymbol{c}$ results as a rational linear
combination of the $n$ most efficient hypergeometric series that satisfy all
stated conditions.

\subsection{Algorithms}

Three algorithms (i) brute force, (ii) MILP solvers and (iii)\textit{\ }%
LLL--Monte Carlo approach, have been implemented by coding PARI--GP\vspace
*{-1pt}\ scripts to compute the solution data matrix (or data table) ${}%
_{b}\boldsymbol{X}_{\boldsymbol{q}}$ of the posed problem\vspace*{-1pt}. The
final multi-valuated formula is obtained by computing $_{b}\boldsymbol{S}%
_{\boldsymbol{q}}$ straightwise from $_{b}\boldsymbol{X}_{\boldsymbol{q}}$ and
simply by inverting this $n\,${\scriptsize x}$\,n$ matrix to get the linear
combination rational factors.

\subsubsection{Brute Force}

In this case the whole lattice Eq.(\ref{14}) $-X_{j}\leq x_{j}\leq
X_{j}$ is scanned for each integer $x_{j},$ $j=1,...n.$ Starting with an empty matrix,
all constraints are evaluated and, if they are satisfied, the feasible
$\boldsymbol{x}$ vector is appended as the $i^{th}$ row $\boldsymbol{x}_{i}$.
Once all nested iteration loops are terminated an $m\,${\scriptsize x}%
$\,n$\ matrix with rows $\boldsymbol{x}_{i},\ i=1,2,...m$ has been collected.
Duplicated rows are removed and the remaining $\boldsymbol{x}_{i}$ are sorted
by objective function value $\varepsilon_{i}=\boldsymbol{x}_{i}\cdot
\boldsymbol{c}$, namely according to the sequence Eq.(\ref{17}). The first top
$n$ rows that give a square matrix of rank $n$ are selected as $_{b}\boldsymbol{X}%
_{\boldsymbol{q}}.$ This method is very fast, simple and robust whenever the
nested iteration deepness $n<8.$ For larger $n $ it hardly works since the
number of iterations increases as $\sim2^{n-1}\Pi_{j=1}^{n}X_{j}$. The brute
force method is coded in PARI function $\mathit{SearchLog}$ with SCIP option
0. (See PARI help\texttt{\
$>$%
?SearchLog}).

\subsubsection{MILP Solvers}

Mixed Integer Linear Programming (MILP) \cite{MILP} \cite{MINLP} is a mathematical NP--hard
optimization method where the objective function and constraints are linear,
but some or all variables are restricted to be integers. In this case only
integer decision variables are used for the linear problem with
Eq.(\ref{16}) excluded. An attempt to exceed the scope of the brute force
method for larger $n$ was performed by implementing algorithms in PARI based on external MILP solvers. 
Several software tools were analyzed, non-commercial and
highly efficient SCIP\texttrademark\ \cite{SCIP}, open-source
HIGHS\texttrademark\ \cite{HIGHS}, commercial MOSEK\texttrademark%
\ \cite{MOSEK} and COPT\texttrademark\ \cite{COPT} among others. SCIP and COPT applications were chosen because they can be easily
controlled inside the PARI-GP script. They are called by the function
$\mathit{SearchLog}$ with SCIP option 1 and 2 respectively. (See script help \texttt{%
$>$%
?SearchLog}).

\bigskip

\hspace*{-15pt}Since MILP solvers cannot compute the whole solutions pool at
once --because just global minima are looked for--, feasible solutions are
obtained for$\ i=1,2,...m$ replacing in Eqs(\ref{13}--\ref{14}) $(\varepsilon
,\boldsymbol{x},x_{j})$ by$\ (\varepsilon_{i},\boldsymbol{x}_{i},x_{ij})\ $and
solving a sequence of ILP extended problems with a new additional linear
constraint given by this objective function lower bound,
\begin{equation}
\varepsilon_{i-1}+\delta\leq\boldsymbol{x}_{i}\cdot\boldsymbol{c}%
\tag{20}\label{20}%
\end{equation}
where $\varepsilon_{0}=0,$ $\varepsilon_{i}=\boldsymbol{x}_{i}\cdot
\boldsymbol{c},$ $i=1,2,...m$ gives the optimal objective sequence pool Eq.(\ref{17} and
$\delta$ is a given absolute pre{\small --}estimated gap ($\delta
=0.5\,${\scriptsize x}$\,10^{-5}$ --an internal SCIP/COPT fixed step
granularity or tolerance limit-- is used) so that this unequality holds,%
\begin{equation}
0<\delta<\min_{1\leq i\leq m}\,(\varepsilon_{i}-\varepsilon_{i-1}%
)\tag{21}\label{21a}%
\end{equation}
Starting with an empty matrix, the ILP solution $\boldsymbol{x}_{i}$ is
appended as a row if the non-linear constraint Eq.(\ref{16}) is satisfied,
otherwise it is rejected. The resulting $\varepsilon_{i}$ are already sorted
so that the $i$--loop simply terminates when the matrix of collected rows gets
rank $n$ and $_{b}\boldsymbol{X}_{\boldsymbol{q}}$ is output. About the
performance, the low acceptance/rejection ratio makes MILP solver method
slightly slower than brute force and the lack of a better precision control
for $\delta$ without a finer granularity (that is an external software
limitation) prevents achieving a larger range for $n.$ The following
algorithm can handle much better these issues.

\subsubsection{LLL--Monte Carlo approach}

A customized Monte Carlo method was implemented to get the pool of best
feasible solutions for this non-linear integer programming problem with larger
$n.$ This is done by using $lindep$\textit{,} the PARI\ GP fast version of the
LLL algorithm \cite{LLL}, that provides a fine tuning control of $d_{s}$, the
internal LLL decimal digits precision (it is specified as an input parameter).
For a sample of size $m,$ $i=1,2,...m$ experiments are performed by setting an
extended random sample vector of logarithms $\boldsymbol{c}_{i}$, defined as
\begin{equation}%
\begin{array}
[c]{c}%
\boldsymbol{c}_{i}=[w_{i},\log p_{1},\log p_{2},...\log p_{n}]%
\genfrac{}{}{0pt}{0}{\genfrac{}{}{0pt}{1}{{}}{{}}}{{}}%
\\
w_{i}=u_{i}\cdot f_{p}\cdot10^{k}%
\genfrac{}{}{0pt}{0}{\genfrac{}{}{0pt}{1}{{}}{{}}}{{}}%
\end{array}
\tag{22}\label{22}%
\end{equation}
where $w_{i}$ is the random component.\ $u_{i}$ is sampled from an uniform
distribution \texttt{U[0,1]}, $f_{p}$ is the machine precision, i.e. the
minimal non--zero positive floating point number and 10$^{k}$ is a scaling
factor with $k\in%
\mathbb{N}
$ a small control integer. The sample values $w_{i}$ are required to be some
orders of magnitude smaller than logarithms while the LLL algorithm must be
able to catch all feasible solution candidates. This integer detection trap is
controlled by means of $k$ and the input precision $d_{s}$. $lindep$ output is
$z_{i},\boldsymbol{x}_{i}$ where $z_{i}$ is an ignored dummy
integer and $\boldsymbol{x}_{i}$ the $i^{th}$ element of the collected
solution candidates vector.
\begin{equation}
\lbrack z_{i},x_{i1},x_{i2},...x_{in}\,]\,=\,lindep(\boldsymbol{c}_{i}%
,d_{s}),\ \ \small{\text{\ }}\ i=1,2,...m\tag{23}\label{23}%
\end{equation}
This algorithm is coded in module $\mathit{SearchLLLog}$ (see PARI help
\texttt{%
$>$%
?SearchLLLog}). For a given $n$ and required input bits $b,$ algorithm starts with a
working global precision about $2b$ and proceeds computing the value of the
internal LLL precision $d_{s}$ (it can be also specified as an input negative
parameter $-\textit{scale}$ in $\mathit{SearchLLLog}$ calls) and $k\!\in%
\mathbb{N}
$, the lowest integer\ in Eq.(\ref{22}) that gives non-trivial LLL integer
detections $lindep(\boldsymbol{c}_{i},d_{s})\neq\lbrack\delta_{0,t}]_{t=0}%
^{n}.$ For a given sample size $m$ (\textit{nmax }in $\mathit{SearchLLLog}$)
nested $i,j,\ell$--loops replacing $d_{s}$\ by $d_{s}+j$ and $k$ by $k+\ell$
with$\ j,\ell=0,1,2,3$ and $i=1,2,...m$ are performed to catch and collect
$\boldsymbol{x}_{i}$ so that, if $m$ is large enough, all relevant
$\boldsymbol{x}_{i}$ are captured in this $4\,${\scriptsize x}$\,4\,${\scriptsize x}$\,m$ polyhedral lattice-net with high probability. When all
loops end, $\boldsymbol{x}_{i}$ duplicates are removed and the remaining
vectors are filtered by non-linear constraint Eq.(\ref{16}). All
$\boldsymbol{x}_{i}$ complying it are sorted by increasing objective function
values $\varepsilon_{i}=\boldsymbol{x}_{i}\cdot\boldsymbol{c}$, Eqs.(\ref{17}%
--\ref{18}). The $n\,${\scriptsize x}$\,n$ solution $_{b}\boldsymbol{X}%
_{\boldsymbol{q}}$ is finally built up as usual with the first top $n$ rows
that give a full rank matrix.

\section{Fast series for single logarithms}

\medskip

In this section, new efficient formulas for logarithms are introduced by
applying the previous integer optimization algorithms. These formulas are
composed by linear combinations of two or more highly efficient hypergeometric
series, resulting in combined identities that are consistently faster and more
effective --sometimes significantly so-- than the rapid single--series
logarithm formulas presented in Part\daurl{:I}{https://arxiv.org/pdf/2506.08245}\aurl{:I}{}. The new formulas span also on a
wider range of logarithms of natural numbers. From now on this convention is
used: $\boldsymbol{\pi}$\texttt{(n)}$=\ $\texttt{[2,3,5,...,p}%
$_{\text{\texttt{n}}}$\texttt{]} is a vector with the complete sequence of the
first $n$ primes; $_{b}\boldsymbol{S}_{\boldsymbol{q}}=$\texttt{[}%
$\mathcal{\,S}_{1},$ $\mathcal{S}_{2},...\mathcal{S}_{n}\,$\texttt{]}%
$^{\text{{\tiny T}}}$ is a column vector representing these $n$ Ramanujan-type
series written as the equivalent form --Eq.(26) Part\daurl{:I}{https://arxiv.org/pdf/2506.08245}\aurl{:I}{}--. (See also Section 10 $\it{Appendix}$).
\begin{equation}
\mathcal{S}_{i}=\frac{1}{\gamma_{i}}\cdot%
{\displaystyle\sum\limits_{k=1}^{\infty}}
\frac{\alpha_{i}\hspace{1pt}k+\beta_{i}}{k(2k-1)}\cdot\left(  \frac
{\mathcal{\nu}_{i}}{\mathcal{\delta}_{i}}\right)  ^{k}\cdot%
\begin{bmatrix}%
\genfrac{}{}{0pt}{0}{{}}{{}}%
1 & \frac{1}{2}%
\genfrac{}{}{0pt}{0}{{}}{{}}%
\\%
\genfrac{}{}{0pt}{0}{{}}{{}}%
\frac{1}{6} & \frac{5}{6}%
\genfrac{}{}{0pt}{0}{{}}{{}}%
\end{bmatrix}
_{k}\hspace{-4pt},\ \ \small{\text{\ }}\ i=1,2,...n\tag{24}\label{24}%
\end{equation}
where all parameters $\alpha_{i},\beta_{i},\gamma_{i},\nu_{i},\delta_{i}\in%
\mathbb{Z}$. $\mathcal{\nu}_{i}\ $and $\mathcal{\delta}_{i}\ne 0$ are
respectively the numerator and denominator of $\rho_{i}$, the rational
convergence rate of the $i^{th}$ series. By abuse of notation $_{b}%
\boldsymbol{S}_{\boldsymbol{q}}$\ is also written as a mapped array%
\begin{equation}
_{b}\boldsymbol{S}_{\boldsymbol{q}}=\left[
\begin{tabular}
[c]{c}%
$\mathcal{S}_{1}$\\
$\mathcal{S}_{2}$\\
$\vdots$\\
$\mathcal{S}_{n}$%
\end{tabular}
\right]  \vspace*{4pt}\leftrightarrow\left[
\begin{tabular}
[c]{ccccc}%
$\mathtt{\alpha}_{1}$ & $\mathtt{\beta}_{1}$ & $\mathtt{\gamma}_{1}$ &
$\mathcal{\nu}_{1}$ & $\mathcal{\delta}_{1}$\\
$\mathtt{\alpha}_{2}$ & $\mathtt{\beta}_{2}$ & $\mathtt{\gamma}_{2}$ &
$\mathcal{\nu}_{2}$ & $\mathcal{\delta}_{2}$\\
$\vdots$ & $\vdots$ & $\vdots$ & $\vdots$ & $\vdots$\\
$\mathtt{\alpha}_{n}$ & $\mathtt{\beta}_{n}$ & $\mathtt{\gamma}_{n}$ &
$\mathcal{\nu}_{n}$ & $\mathcal{\delta}_{n}$%
\end{tabular}
\right] \tag{25}\label{25}%
\end{equation}
Note that $\alpha_{i},\beta_{i},\gamma_{i},\nu_{i},\delta_{i}$ are usually large, so $_{b}\boldsymbol{S}%
_{\boldsymbol{q}}$\ formulas are given explicitly just for $b=64\,$%
\texttt{-bit}. Formulas with higher input precision, $b>64$, are given only through the
solution matrix $_{b}\boldsymbol{X}_{\boldsymbol{q}}$ that just contains small
integers, in this case the final column vectors of series and logarithms are
obtained applying Eqs.(\ref{1}--\ref{2}) or Eq.(\ref{9}) and Eq.(\ref{19}).
With this setup some new multi-valuation formulas giving arbitrary precision output, are as follows,

\subsection{log(2) and log(3)}

To get \texttt{64-bit} formulas in this case, the application of coded
functions \texttt{SearchLog([2,3],[],64,0.6)}\ or
\texttt{SearchLLLog([2,3],64,1/4,128)} gives a total binary splitting cost
\texttt{C}$_{\text{s}}\ $\texttt{=}\ \texttt{0.759456}, significatively better
than costs of single series Eqs.(18--19) Part\daurl{:I}{https://arxiv.org/pdf/2506.08245}\aurl{:I}{}, the fastest formulas
for these constants known so far. Best formulas are now,
\begin{equation}%
\begin{tabular}
[c]{rrrrrrr}%
${}_{\text{\texttt{64}}}\boldsymbol{X}_{\boldsymbol{\pi}\text{\texttt{(2)}}} $
& $=$ & $\left[
\begin{tabular}
[c]{rr}%
\texttt{8} & \texttt{-5}\\
\texttt{3} & \texttt{-2}%
\end{tabular}
\right]  $ &  & $_{\text{\texttt{64}}}\boldsymbol{X}_{\boldsymbol{\pi
}\text{\texttt{(2)}}}^{-1}$ & $=$ & $\left[
\begin{tabular}
[c]{rr}%
\texttt{2} & \texttt{-5}\\
\texttt{3} & \texttt{-8}%
\end{tabular}
\right]  $%
\end{tabular}
{}\tag{26}\label{26}%
\end{equation}
$\mathtt{\ }$%
\[%
\genfrac{}{}{0pt}{0}{_{\text{\texttt{64}}}\boldsymbol{S}_{\boldsymbol{\pi
}\text{\texttt{(2)}}}\genfrac{}{}{0pt}{0}{\genfrac{}{}{0pt}{1}{{}}{{}}}{{}%
}}{\left[
\begin{tabular}
[c]{rrrrr}%
\texttt{278133806980282} & \texttt{-46355624995421} & \texttt{742586} &
\texttt{4826809} & \texttt{104068027861696512}\\
\texttt{12705086} & \texttt{-2117503} & \texttt{-2} & \texttt{1} &
\texttt{161803008}%
\end{tabular}
\right]  }%
\vspace*{5pt}%
\]
finally,%
\begin{equation}%
\begin{tabular}
[c]{c}%
\texttt{log(2)\ =\ 2\thinspace S}$_{\text{\texttt{1}}}\ $%
\texttt{-\ 5\thinspace S}$_{\text{\texttt{2}}}%
\genfrac{}{}{0pt}{0}{\genfrac{}{}{0pt}{1}{{}}{{}}}{{}}%
$\\
\texttt{log(3)\ =\ 3\thinspace S}$_{\text{\texttt{1}}}\ $%
\texttt{-\ 8\thinspace S}$_{\text{\texttt{2}}}%
\genfrac{}{}{0pt}{0}{\genfrac{}{}{0pt}{1}{{}}{{}}}{{}}%
$%
\end{tabular}
\tag{27}\label{27}%
\end{equation}

\subsection{log(2), log(5) and log(10)}

In this case for \texttt{64-bit} formulas the application of functions
\texttt{SearchLog([2,5],[],64,0.6)}\ or \texttt{SearchLLLog([2,5],64,1/4,128)}
gives a total binary splitting cost \texttt{C}$_{\text{s}}\ $%
\texttt{=\ 0.811450}, significatively better than single series Eq.(20)
Part\daurl{:I}{https://arxiv.org/pdf/2506.08245}\aurl{:I}{} cost. Best formulas are%
\begin{equation}%
\begin{tabular}
[c]{rrrrrrr}%
${}_{\text{\texttt{64}}}\boldsymbol{X}_{\text{\texttt{[2,5]}}}$ & $=$ &
$\left[
\begin{tabular}
[c]{rr}%
\texttt{-7} & \texttt{3}\\
\texttt{-2} & \texttt{1}%
\end{tabular}
\right]  $ &  & $_{\text{\texttt{64}}}\boldsymbol{X}_{\text{\texttt{[2,5]}}%
}^{-1}$ & $=$ & $\left[
\begin{tabular}
[c]{rr}%
\texttt{-1} & \texttt{3}\\
\texttt{-2} & \texttt{7}%
\end{tabular}
\right]  $%
\end{tabular}
{}\tag{28}\label{28}%
\end{equation}%
\begin{equation}%
\genfrac{}{}{0pt}{0}{_{\text{\texttt{64}}}\boldsymbol{S}_{\text{\texttt{[2,5]}%
}}\genfrac{}{}{0pt}{0}{\genfrac{}{}{0pt}{1}{{}}{{}}}{{}}}{\left[
\begin{tabular}
[c]{rrrrr}%
\texttt{9327029143014} & \texttt{-1554504843507} & \texttt{-486} & \texttt{27}
& \texttt{65545216000000}\\
\texttt{520542} & \texttt{-86751} & \texttt{2} & \texttt{1} & \texttt{3499200}%
\end{tabular}
\right]  }%
\vspace*{5pt}\tag{29}\label{29}%
\end{equation}
and%
\begin{equation}%
\begin{tabular}
[c]{l}%
\texttt{log(10)\ =\ -3\thinspace S}$_{\text{\texttt{1}}}\ $%
\texttt{+\ 10\thinspace S}$_{\text{\texttt{2}}}%
\genfrac{}{}{0pt}{0}{\genfrac{}{}{0pt}{1}{{}}{{}}}{{}}%
$\\
\texttt{log(5)\ \ =\ -2\thinspace S}$_{\text{\texttt{1}}}\ $\texttt{+
\ 7\thinspace S}$_{\text{\texttt{2}}}%
\genfrac{}{}{0pt}{0}{\genfrac{}{}{0pt}{1}{{}}{{}}}{{}}%
$\\
\texttt{log(2)\ \ =\ \ -\thinspace S}$_{\text{\texttt{1}}}\ $\texttt{+
\ 3\thinspace S}$_{\text{\texttt{2}}}%
\genfrac{}{}{0pt}{0}{\genfrac{}{}{0pt}{1}{{}}{{}}}{{}}%
$%
\end{tabular}
\tag{30}\label{30}%
\end{equation}

\subsection{log(7)}

\texttt{SearchLog([2,3,7],,64,0.6)}\ or
\texttt{SearchLLLog([2,3,7],64,1/4,128)\ }give \texttt{64-bit} formulas to
compute \texttt{log(7)} with total binary splitting cost \texttt{C}%
$_{\text{s}}\ $\texttt{=\ 0.909610}, in this case%
\begin{equation}%
\begin{tabular}
[c]{rrr}%
${}_{\text{\texttt{64}}}\boldsymbol{X}_{\text{\texttt{[2,3,7]}}}=\left[
\begin{tabular}
[c]{rrr}%
\texttt{-6} & \texttt{2} & \texttt{1}\\
\texttt{-4} & \texttt{-1} & \texttt{2}\\
\texttt{-5} & \texttt{5} & \texttt{-1}%
\end{tabular}
\right]  $ &  & $_{\text{\texttt{64}}}\boldsymbol{X}_{\text{\texttt{[2,3,7]}}%
}^{-1}=\left[
\begin{tabular}
[c]{rrr}%
\texttt{-9} & \texttt{7} & \texttt{5}\\
\texttt{-14} & \texttt{11} & \texttt{8}\\
\texttt{-25} & \texttt{20} & \texttt{14}%
\end{tabular}
\right]  $%
\end{tabular}
{}\tag{31}\label{31}%
\end{equation}%
\[%
\genfrac{}{}{0pt}{0}{_{\text{\texttt{64}}}\boldsymbol{S}%
_{\text{\texttt{[2,3,7]}}}\genfrac{}{}{0pt}{0}{\genfrac{}{}{0pt}{1}{{}}{{}%
}}{{}}}{\left[
\begin{tabular}
[c]{rrrrr}%
\texttt{297314599426} & \texttt{-49552433153} & \texttt{-2} & \texttt{1} &
\texttt{28318630330368}\\
\texttt{77272372606} & \texttt{-12878728703} & \texttt{2} & \texttt{1} &
\texttt{5621365951488}\\
\texttt{199355237389946} & \texttt{-33225832325053} & \texttt{4952198} &
\texttt{47045881} & \texttt{69785645582757888}%
\end{tabular}
\right]  }%
\vspace*{8pt}%
\]
so that,\vspace*{4pt}%

\begin{equation}%
\begin{tabular}
[c]{c}%
\texttt{log(7)\ =\ -25\thinspace S}$_{\text{\texttt{1}}}\ $%
\texttt{+\ 20\thinspace S}$_{\text{\texttt{2}}}\ $\texttt{+\ 14\thinspace
S}$_{\text{\texttt{3}}}$%
\end{tabular}
\vspace*{4pt}\tag{32}\label{32}%
\end{equation}

\subsection{log(11)}

\texttt{SearchLog([2,3,11],,64,0.6)}\ or
\texttt{SearchLLLog([2,3,11],64,,128)\ }give \texttt{64-bit} formulas to
compute \texttt{log(11)} with total binary splitting cost \texttt{C}%
$_{\text{s}}\ $\texttt{=\ 0.838056}, in this case
\begin{equation}%
\begin{tabular}
[c]{rrr}%
${}_{\text{\texttt{64}}}\boldsymbol{X}_{\text{\texttt{[2,3,11]}}}=\left[
\begin{tabular}
[c]{rrr}%
\texttt{-1} & \texttt{5} & \texttt{-2}\\
\texttt{-5} & \texttt{1} & \texttt{1}\\
\texttt{-8} & \texttt{5} & \texttt{0}%
\end{tabular}
\right]  $ &  & $_{\text{\texttt{64}}}\boldsymbol{X}_{\text{\texttt{[2,3,11]}%
}}^{-1}=\left[
\begin{tabular}
[c]{rrr}%
\texttt{5} & \texttt{10} & \texttt{-7}\\
\texttt{8} & \texttt{16} & \texttt{-11}\\
\texttt{17} & \texttt{35} & \texttt{-24}%
\end{tabular}
\right]  $%
\end{tabular}
{}\tag{33}\label{33}%
\end{equation}%
\[%
\genfrac{}{}{0pt}{0}{_{\text{\texttt{64}}}\boldsymbol{S}%
_{\text{\texttt{[2,3,11]}}}\genfrac{}{}{0pt}{0}{\genfrac{}{}{0pt}{1}{{}}{{}%
}}{{}}}{\left[
\begin{tabular}
[c]{rrrrr}%
\texttt{241517233468190} & \texttt{-40252872244375} & \texttt{2} & \texttt{1}
& \texttt{87851769180634800}\\
\texttt{10438496510} & \texttt{-1739749375} & \texttt{2} & \texttt{1} &
\texttt{508836556800}\\
\texttt{278133806980282} & \texttt{-46355624995421} & \texttt{-742586} &
\texttt{4826809} & \texttt{104068027861696512}%
\end{tabular}
\right]  }%
\vspace*{8pt}%
\]%
\begin{equation}%
\begin{tabular}
[c]{c}%
\texttt{log(11)\ =\ 17\thinspace S}$_{\text{\texttt{1}}}\ $%
\texttt{+\ 35\thinspace S}$_{\text{\texttt{2}}}\ $\texttt{-\ 24\thinspace
S}$_{\text{\texttt{3}}}$%
\end{tabular}
\vspace*{4pt}\tag{34}\label{34}%
\end{equation}

\subsection{log(13)}

\texttt{SearchLog([2,3,13],,64,0.6)}\ or
\texttt{SearchLLLog([2,3,13],64,,128)\ }give \texttt{64-bit} formulas to
compute \texttt{log(13)} with total binary splitting cost \texttt{C}%
$_{\text{s}}\ $\texttt{=\ 0.966545}, so that
\begin{equation}%
\begin{tabular}
[c]{rrr}%
${}_{\text{\texttt{64}}}\boldsymbol{X}_{\text{\texttt{[2,3,13]}}}=\left[
\begin{tabular}
[c]{rrr}%
\texttt{-1} & \texttt{3} & \texttt{-1}\\
\texttt{-1} & \texttt{-4} & \texttt{2}\\
\texttt{-8} & \texttt{5} & \texttt{0}%
\end{tabular}
\right]  $ &  & $_{\text{\texttt{64}}}\boldsymbol{X}_{\text{\texttt{[2,3,13]}%
}}^{-1}=\left[
\begin{tabular}
[c]{rrr}%
\texttt{10} & \texttt{5} & \texttt{-2}\\
\texttt{16} & \texttt{8} & \texttt{-3}\\
\texttt{37} & \texttt{19} & \texttt{-7}%
\end{tabular}
\right]  $%
\end{tabular}
{}\tag{35}\label{35}%
\end{equation}%
\[%
\genfrac{}{}{0pt}{0}{_{\text{\texttt{64}}}\boldsymbol{S}%
_{\text{\texttt{[2,3,13]}}}\genfrac{}{}{0pt}{0}{\genfrac{}{}{0pt}{1}{{}}{{}%
}}{{}}}{\left[
\begin{tabular}
[c]{rrrrr}%
\texttt{3761526494} & \texttt{-626921047} & \texttt{2} & \texttt{1} &
\texttt{149502935088}\\
\texttt{35732110926898} & \texttt{\ -5955351291329} & \texttt{33614} &
\texttt{117649} & \texttt{8869174125759792}\\
\texttt{278133806980282} & \texttt{-46355624995421} & \texttt{-742586} &
\texttt{4826809} & \texttt{104068027861696512}%
\end{tabular}
\right]  }%
\vspace*{8pt}%
\]%
\begin{equation}%
\begin{tabular}
[c]{lcl}%
\texttt{log(13)} & \texttt{=} & \texttt{37\thinspace S}$_{\text{\texttt{1}}%
}\ $\texttt{+\ 19\thinspace S}$_{\text{\texttt{2}}}\ $\texttt{-\ 7\thinspace
S}$_{\text{\texttt{3}}}$%
\end{tabular}
\vspace*{0pt}\tag{36}\label{36}%
\end{equation}

\subsection{log(3), log(5), log(10) and log(15)}

\texttt{64-bit} formulas for these logarithms with total binary splitting cost
\texttt{C}$_{\text{s}}\ $\texttt{=\ 0.819035} are obtained by
\texttt{SearchLog([2,3,5],,64,0.6)} or
\texttt{SearchLLLog([2,3,5],64,,128)\vspace*{4pt}} as linear combinations of three series,
\begin{equation}%
\begin{tabular}
[c]{rrr}%
${}_{\text{\texttt{64}}}\boldsymbol{X}_{\boldsymbol{\pi}\text{\texttt{(3)}}%
}=\left[
\begin{tabular}
[c]{rrr}%
\texttt{-4} & \texttt{4} & \texttt{1}\\
\texttt{-7} & \texttt{0} & \texttt{3}\\
\texttt{-1} & \texttt{5} & \texttt{-3}%
\end{tabular}
\right]  $ &  & $_{\text{\texttt{64}}}\boldsymbol{X}_{\boldsymbol{\pi
}\text{\texttt{(3)}}}^{-1}=\left[
\begin{tabular}
[c]{rrr}%
\texttt{15} & \texttt{-7} & \texttt{-12}\\
\texttt{24} & \texttt{-11} & \texttt{-19}\\
\texttt{35} & \texttt{-16} & \texttt{-28}%
\end{tabular}
\right]  $%
\end{tabular}
{}\tag{37}\label{37}%
\end{equation}%
\[%
\genfrac{}{}{0pt}{0}{_{\text{\texttt{64}}}\boldsymbol{S}_{\boldsymbol{\pi
}\text{\texttt{(3)}}}\genfrac{}{}{0pt}{0}{\genfrac{}{}{0pt}{1}{{}}{{}}}{{}%
}}{\left[
\begin{tabular}
[c]{rrrrr}%
\texttt{973517952638} & \texttt{-162252991999} & \texttt{2} & \texttt{1} &
\texttt{117550781107200}\\
\texttt{9327029143014} & \texttt{\ -1554504843507} & \texttt{-486} &
\texttt{27} & \texttt{65545216000000}\\
\texttt{262018021085614} & \texttt{-43669669391807} & \texttt{-33614} &
\texttt{117649} & \texttt{96874652706750000}%
\end{tabular}
\right]  }%
\vspace*{8pt}%
\]%
\begin{equation}%
\begin{tabular}
[c]{lcl}%
\texttt{log(3)} & \texttt{=} & \texttt{24\thinspace S}$_{\text{\texttt{1}}}%
\ $\texttt{-\ 11\thinspace S}$_{\text{\texttt{2}}}\ $\texttt{-\ 19\thinspace
S}$_{\text{\texttt{3}}}%
\genfrac{}{}{0pt}{0}{\genfrac{}{}{0pt}{1}{{}}{{}}}{{}}%
$\\
\texttt{log(5)} & \texttt{=} & \texttt{35\thinspace S}$_{\text{\texttt{1}}}%
\ $\texttt{-\ 16\thinspace S}$_{\text{\texttt{2}}}\ $\texttt{-\ 28\thinspace
S}$_{\text{\texttt{3}}}%
\genfrac{}{}{0pt}{0}{\genfrac{}{}{0pt}{1}{{}}{{}}}{{}}%
$\\
\texttt{log(10)} & \texttt{=} & \texttt{50\thinspace S}$_{\text{\texttt{1}}%
}\ $\texttt{-\ 23\thinspace S}$_{\text{\texttt{2}}}\ $\texttt{-\ 40\thinspace
S}$_{\text{\texttt{3}}}%
\genfrac{}{}{0pt}{0}{\genfrac{}{}{0pt}{1}{{}}{{}}}{{}}%
$\\
\texttt{log(15)} & \texttt{=} & \texttt{59\thinspace S}$_{\text{\texttt{1}}%
}\ $\texttt{-\ 27\thinspace S}$_{\text{\texttt{2}}}\ $\texttt{-\ 47\thinspace
S}$_{\text{\texttt{3}}}%
\genfrac{}{}{0pt}{0}{\genfrac{}{}{0pt}{1}{{}}{{}}}{{}}%
$%
\end{tabular}
\vspace*{4pt}\tag{38}\label{38}%
\end{equation}

\subsection{log(17)}

\texttt{SearchLog([2,3,17],,64,0.6)}\ or
\texttt{SearchLLLog([2,3,17],64,,128)}\ give \texttt{64-bit} formulas to
compute \texttt{log(17)} with total binary splitting cost \texttt{C}%
$_{\text{s}}\ $\texttt{= 0.880169}$\ $, so that\vspace*{-6pt}
\begin{equation}%
\begin{tabular}
[c]{rrr}%
${}_{\text{\texttt{64}}}\boldsymbol{X}_{\text{\texttt{[2,3,17]}}}=\left[
\begin{tabular}
[c]{rrr}%
\texttt{-5} & \texttt{-2} & \texttt{2}\\
\texttt{-8} & \texttt{5} & \texttt{0}\\
\texttt{-1} & \texttt{-2} & \texttt{1}%
\end{tabular}
\right]  $ &  & $_{\text{\texttt{64}}}\boldsymbol{X}_{\text{\texttt{[2,3,17]}%
}}^{-1}=\left[
\begin{tabular}
[c]{rrr}%
\texttt{5} & \texttt{-2} & \texttt{-10}\\
\texttt{8} & \texttt{-3} & \texttt{-16}\\
\texttt{21} & \texttt{-8} & \texttt{-41}%
\end{tabular}
\right]  $%
\end{tabular}
{}\tag{39}\label{39}%
\end{equation}
$\vspace*{-5pt}$%
\[%
\genfrac{}{}{0pt}{0}{_{\text{\texttt{64}}}\boldsymbol{S}%
_{\text{\texttt{[2,3,17]}}}\genfrac{}{}{0pt}{0}{\genfrac{}{}{0pt}{1}{{}}{{}%
}}{{}}}{\left[
\begin{tabular}
[c]{rrrrr}%
\texttt{575598165481726} & \texttt{-95933027579903} & \texttt{2} & \texttt{1}
& \texttt{249089856719597568}\\
\texttt{278133806980282} & \texttt{-46355624995421} & \texttt{-742586} &
\texttt{4826809} & \texttt{104068027861696512}\\
\texttt{472053890} & \texttt{-78675625} & \texttt{-2} & \texttt{1} &
\texttt{12388042800}%
\end{tabular}
\right]  }%
\]%
\begin{equation}%
\begin{tabular}
[c]{c}%
\texttt{log(17)\ =\ 21\thinspace S}$_{\text{\texttt{1}}}\ $%
\texttt{-\ 8\thinspace S}$_{\text{\texttt{2}}}\ $\texttt{-\ 41\thinspace
S}$_{\text{\texttt{3}}}$%
\end{tabular}
\vspace*{3pt}\tag{40}\label{40}%
\end{equation}

\subsection{log(19)}

\texttt{SearchLog([2,3,5,19],,64,.6),}%
\ \texttt{SearchLLLog([2,3,5,19],64,,128)}\ give \texttt{64-bit} formulas to
compute \texttt{log(19)} with total binary splitting cost \texttt{C}%
$_{\text{s}}\ $\texttt{=\ 0.971132}, in this case
\begin{equation}%
\begin{tabular}
[c]{ccc}%
$_{\text{\texttt{64}}}\boldsymbol{X}_{\text{\texttt{[2,3,5,19]}}}$ &  &
$_{\text{\texttt{64}}}\boldsymbol{X}_{\text{\texttt{[2,3,5,19]}}}^{-1}$\\
& \multicolumn{1}{r}{} & \\
\multicolumn{1}{r}{$\left[
\begin{tabular}
[c]{rrrr}%
\texttt{-3} & \texttt{-2} & \texttt{-1} & \texttt{2}\\
\texttt{-5} & \texttt{-1} & \texttt{1} & \texttt{1}\\
\texttt{-4} & \texttt{4} & \texttt{-1} & \texttt{0}\\
\texttt{-1} & \texttt{5} & \texttt{-3} & \texttt{0}%
\end{tabular}
\right]  $} & \multicolumn{1}{r}{} & \multicolumn{1}{r}{$\left[
\begin{tabular}
[c]{rrrr}%
\texttt{7} & \texttt{-14} & \texttt{15} & \texttt{-12}\\
\texttt{11} & \texttt{-22} & \texttt{24} & \texttt{-19}\\
\texttt{16} & \texttt{-32} & \texttt{35} & \texttt{-28}\\
\texttt{30} & \texttt{-59} & \texttt{64} & \texttt{-51}%
\end{tabular}
\right]  $}%
\end{tabular}
{}\tag{41}\label{41}%
\end{equation}
$\vspace*{-1pt}$%
\[%
\genfrac{}{}{0pt}{0}{_{\text{\texttt{64}}}\boldsymbol{S}%
_{\text{\texttt{[2,3,5,19]}}}\genfrac{}{}{0pt}{0}{\genfrac{}{}{0pt}{1}{{}}{{}%
}}{{}}}{\left[
\begin{tabular}
[c]{rrrrr}%
\texttt{1753547120930878} & \texttt{-292257853487999} & \texttt{2} &
\texttt{1} & \texttt{948229997617324800}\\
\texttt{2287649600258} & \texttt{-381274933249} & \texttt{-2} & \texttt{1} &
\texttt{327702810931200}\\
\texttt{973517952638} & \texttt{-162252991999} & \texttt{2} & \texttt{1} &
\texttt{117550781107200}\\
\texttt{262018021085614} & \texttt{-43669669391807} & \texttt{-33614} &
\texttt{117649} & \texttt{96874652706750000}%
\end{tabular}
\right]  }%
\]\vspace*{4pt}%
\begin{equation}%
\begin{tabular}
[c]{c}%
\texttt{log(19)\ =\ 30\thinspace S}$_{\text{\texttt{1}}}\ $%
\texttt{-\ 59\thinspace S}$_{\text{\texttt{2}}}$\texttt{\ +\ 64\thinspace
S}$_{\text{\texttt{3}}}\ $\texttt{-\ 51\thinspace S}$_{\text{\texttt{4}}}$%
\end{tabular}
\vspace*{4pt}\tag{42}\label{42}%
\end{equation}

\subsection{log($n$)}

Function \texttt{SearchLog([2,3,5,n],,b,tol,1/8) }(or a variation) can be
applied to get fast\ \texttt{b-bit} formulas to compute \texttt{log(n)\ }for
larger $n$ delivering low total binary splitting cost \texttt{C}$_{\text{s}}$.
In this case $n$ must have some factor different from \texttt{[2,3,5]}
(otherwise some of $_{\text{\texttt{b}}}\boldsymbol{X}_{\text{\texttt{[2,3,5]}%
}},\ _{\text{\texttt{b}}}\boldsymbol{X}_{\text{\texttt{[2,3]}}}$
or$\ _{\text{\texttt{b}}}\boldsymbol{X}_{\text{\texttt{[2,5]}}}$ is applied).
Using \texttt{b = 192} a wide range of $n$ is covered, up to hundred thousands
or more (\texttt{b} must be larger to go further).

\subsubsection{Example 3}

Consider \texttt{n = prime(1048576)\ = 16290047}. This solution matrix is
found by \texttt{SearchLog([2,3,5,16290047],,192,0.006,1/4)} giving
\texttt{C}$_{\text{s}}\ $\texttt{=\ 0.723618,\vspace*{1pt}}%
\begin{equation}
_{\text{\texttt{192}}}\boldsymbol{X}_{\text{\texttt{[2,3,5,16290047]}}%
}=\left[
\begin{tabular}
[c]{rrrr}%
\texttt{0} & \texttt{-1} & \texttt{11} & \texttt{-1}\\
\texttt{-15} & \texttt{8} & \texttt{1} & \texttt{0}\\
\texttt{-8} & \texttt{-13} & \texttt{2} & \texttt{1}\\
\texttt{-9} & \texttt{13} & \texttt{-5} & \texttt{0}%
\end{tabular}
\right] \tag{43}\label{43}%
\end{equation}
\texttt{\vspace*{1pt}}%
\[%
\begin{tabular}
[c]{c}%
\texttt{log(16290047)\ =\ -\ 1270\thinspace S}$_{\text{\texttt{1}}}%
\ $\texttt{+\ 2373\thinspace S}$_{\text{\texttt{2}}}$%
\texttt{\ -\ 1269\thinspace S}$_{\text{\texttt{3}}}\ $%
\texttt{-\ 2827\thinspace S}$_{\text{\texttt{4}}}$%
\end{tabular}
\vspace*{2pt}%
\]
so that formulas are defined by 16 small integers that are mapped to $n,$
generating very low storage tables.\ \texttt{log(n)} is computed as a linear
combination of 4 ultra fast series {\small --}Eq.(\ref{19}) using either
Eqs.(\ref{1}--\ref{2}) or Eq.(\ref{9}){\small --} that are easily built from
these stored data tables.

\section{Single logarithm formulas testing}

\bigskip
\begin{center}%
\begin{tabular}
[c]{|c|c|c|c|c|}\hline\hline
\texttt{Constant} & $%
\genfrac{}{}{0pt}{}{\text{\texttt{Formula}}}{\text{\texttt{Reference}}}%
$ & \texttt{y-cruncher} & $%
\genfrac{}{}{0pt}{}{\text{\texttt{Machin-like formula}}}{\text{\texttt{x}}%
_{i}\text{\texttt{\ in acoth(x}}_{i}\mathtt{)}}%
$ & $\texttt{alias} $\\\hline\hline
\texttt{log(2)} & \multicolumn{1}{|r|}{\texttt{Part I, Eq.(18)}} &
\texttt{-algorithm:0} & \texttt{-} & \texttt{-alg0}\\
\texttt{log(2)} & \multicolumn{1}{|r|}{\texttt{Part I, Eq.(21)}} &
\texttt{-algorithm:1} & \texttt{-} & \texttt{-alg1}\\
\texttt{log(2)} & \multicolumn{1}{|r|}{\texttt{Eq.(\ref{27})}} & \texttt{-} &
\texttt{-} & \texttt{2025A}\\
\texttt{log(2)} & \multicolumn{1}{|r|}{\texttt{Eq.(\ref{30})}} & \texttt{-} &
\texttt{-} & \texttt{2025B}\\
\hline
\texttt{log(3)} & \multicolumn{1}{|r|}{\texttt{Part I, Eq.(19)}} &
\texttt{-algorithm:0} & \texttt{-} & \texttt{-alg0}\\
\texttt{log(3)} & \multicolumn{1}{|r|}{\texttt{Part I, Eq.(24)}} &
\texttt{-algorithm:1} & \texttt{-} & \texttt{-alg1}\\
\texttt{log(3)} & \multicolumn{1}{|r|}{\texttt{Eq.(\ref{27})}} & \texttt{-} &
\texttt{-} & \texttt{2025A}\\
\texttt{log(3)} & \multicolumn{1}{|r|}{\texttt{Eq.(\ref{38})}} & \texttt{-} &
\texttt{-} & \texttt{2025B}\\
\hline
\texttt{log(5)} & \multicolumn{1}{|r|}{\texttt{Part I, Eq.(20)}} &
\texttt{-algorithm:0} & \texttt{-} & \texttt{-alg0}\\
\texttt{log(5)} & \multicolumn{1}{|r|}{\texttt{\cite{NumberW}}} &
\texttt{-algorithm:1} & \texttt{251,449,4801,8749} & \texttt{-alg1}\\
\texttt{log(5)} & \multicolumn{1}{|r|}{\texttt{Eq.(\ref{30})}} & \texttt{-} &
\texttt{-} & \texttt{2025A}\\
\texttt{log(5)} & \multicolumn{1}{|r|}{\texttt{Eq.(\ref{38})}} & \texttt{-} &
\texttt{-} & \texttt{2025B}\\
\hline
\texttt{log(7)} & \multicolumn{1}{|r|}{\texttt{\cite{NumberX}}} &
\texttt{-algorithm:0} & \texttt{251,449,4801,8749} & \texttt{-alg0}\\
\texttt{log(7)} & \multicolumn{1}{|r|}{\texttt{\cite{NumberX}}} &
\texttt{-algorithm:1} & \texttt{99,449,4801,8749} & \texttt{-alg1}\\
\texttt{log(7)} & \multicolumn{1}{|r|}{\texttt{Eq.(\ref{32})}} & \texttt{-} &
\texttt{-} & \texttt{2025A}\\
\texttt{log(7)} & \multicolumn{1}{|r|}{$\text{\texttt{primes} \texttt{=}
\texttt{[2,5,7]}}$} & \texttt{-} & \texttt{-} & \texttt{2025B}\\
\hline
\texttt{log(10)} & \multicolumn{1}{|r|}{\texttt{\cite{NumberW}}} &
\texttt{-algorithm:0} & \multicolumn{1}{|r|}{\texttt{251,449,4801,8749}} &
\texttt{alg0}\\
\texttt{log(10)} & \multicolumn{1}{|r|}{\texttt{\cite{NumberW}}} &
\texttt{-algorithm:1} & \multicolumn{1}{|r|}{\texttt{99,449,4801,8749}} &
\texttt{-alg1}\\
\texttt{log(10)} & \multicolumn{1}{|r|}{\texttt{Eq.(\ref{30})}} & \texttt{-} &
\texttt{-} & \texttt{2025A}\\
\texttt{log(10)} & \multicolumn{1}{|r|}{\texttt{Eq.(\ref{38})}} & \texttt{-} &
\texttt{-} & \texttt{2025B}\\
\hline
\texttt{log(11)} & \multicolumn{1}{|r|}{\texttt{\cite{NumberX}}} &
\texttt{-algorithm:0} & \multicolumn{1}{|r|}{\texttt{65,485,769,19601}} &
\texttt{-alg0}\\
\texttt{log(11)} & \multicolumn{1}{|r|}{\texttt{\cite{NumberX}}} &
\texttt{-algorithm:1} & \multicolumn{1}{|c|}{\texttt{5 terms}} &
\texttt{-alg1}\\
\texttt{log(11)} & \multicolumn{1}{|r|}{\texttt{Eq.(\ref{34})}} & \texttt{-} &
\texttt{-} & \texttt{2025A}\\
\texttt{log(11)} & $\text{\texttt{primes} \texttt{=} \texttt{[2,5,11]}}$ &
\texttt{-} & \texttt{-} & \texttt{2025B}\\
\hline
\multicolumn{5}{|c|}{\vspace*{-8pt}}\\
\multicolumn{5}{|c|}{\texttt{Table I. Reference Formulas for Testings}}\\\hline
\end{tabular}

\texttt{\vspace*{-4pt}}
\end{center}
\texttt{NOTES:}

{\small -- algorithm:0, \texttt{[-alg0]} y-cruncher's primary native algorithm.}

{\small -- algorithm:1, \texttt{[-alg1]} y-cruncher's secondary native algorithm for digits
verification. }

{\small -- Machin-like, integer arguments of acoth functions in Machin-type
formulas. (Table 1,\thinspace\cite{FREDJ})}

{\small -- primes = }$[p_{1},p_{2},p_{3}]${\small , optimal identities got
by\vspace*{8pt}} \texttt{SearchLog([p}$_{\text{\texttt{1}}}$\texttt{,p}%
$_{\text{\texttt{2}}}$\texttt{,p}$_{\text{\texttt{3}}}$%
\texttt{],,64,0.6,1/4)\vspace*{2pt}}

Custom configuration files for \texttt{64-bit} formulas Eqs.(\ref{28}%
--\ref{42}) were written for y-cruncher software \cite{YCRUNCHER} to compute
\texttt{1B}, \texttt{10B}, \texttt{100B} and \texttt{340B} $(\texttt{1B} = 10^9\  \text{--billion--})$ decimal digits and test their performances against the currently
fastest native logarithm algorithms embedded on this platform. 
We used a \texttt{Dell Precision 7960 Intel Xeon W9-3595X CPU} with 60 cores, 60 threads and \texttt{2TB RAM} running on MS Windows\textsuperscript{TM} 11 Pro for Workstations, 25H2 OS (this setup and source data for figures and graphs in next pages were kindly provided by Dmitriy Grigoryev). These optimal integer programming solutions
give the currently known fastest primary formulas for multiprecision computing
of logarithms of small integers using the binary splitting algorithm. From
these testings the speed improvement of the new identities is small but increasing for $\log(2)$ at a huge number of digits (greater than 100B) however it is significant in the whole range for $\log(3)$, $\log(5)$, $\log(7)$,
$\log(10)$, $\log(11)$ (also for logarithms of higher integers not
displayed in Table I above). The results of these tests are shown below. \pagebreak \newline \newline \\
\\ \\ \\
\begin{figure}[h] 
\raggedright
\includegraphics[viewport=0in 0in 11.552200in 4.052500in,
height=1.77534in,
width=5.01552in
]%
{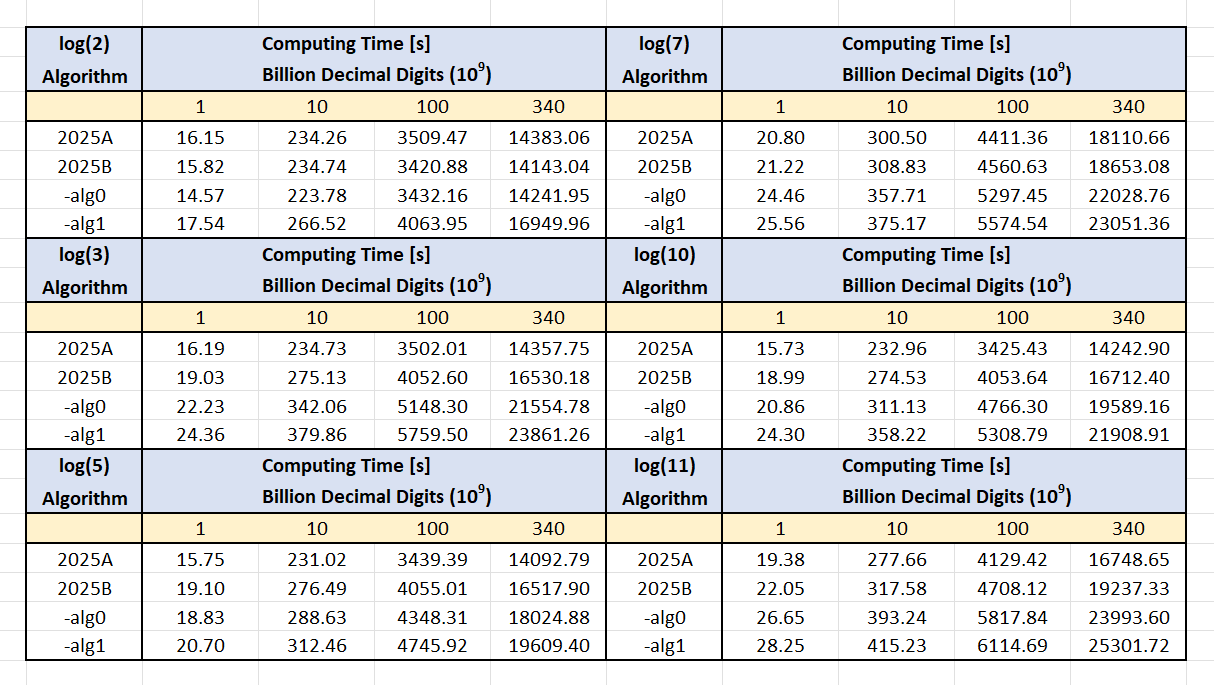}%
\vspace{-6pt}
\caption{Computing Time in \texttt{[s]} for \texttt{log} algorithms.}
\end{figure}
\newline \newline \\
\\ \\ \newline \newline 
\begin{figure}[h] 
	\raggedright
	\includegraphics[viewport=0in 0in 11.552200in 4.052500in,
	height=1.77534in,
	width=5.01552in
	]%
	{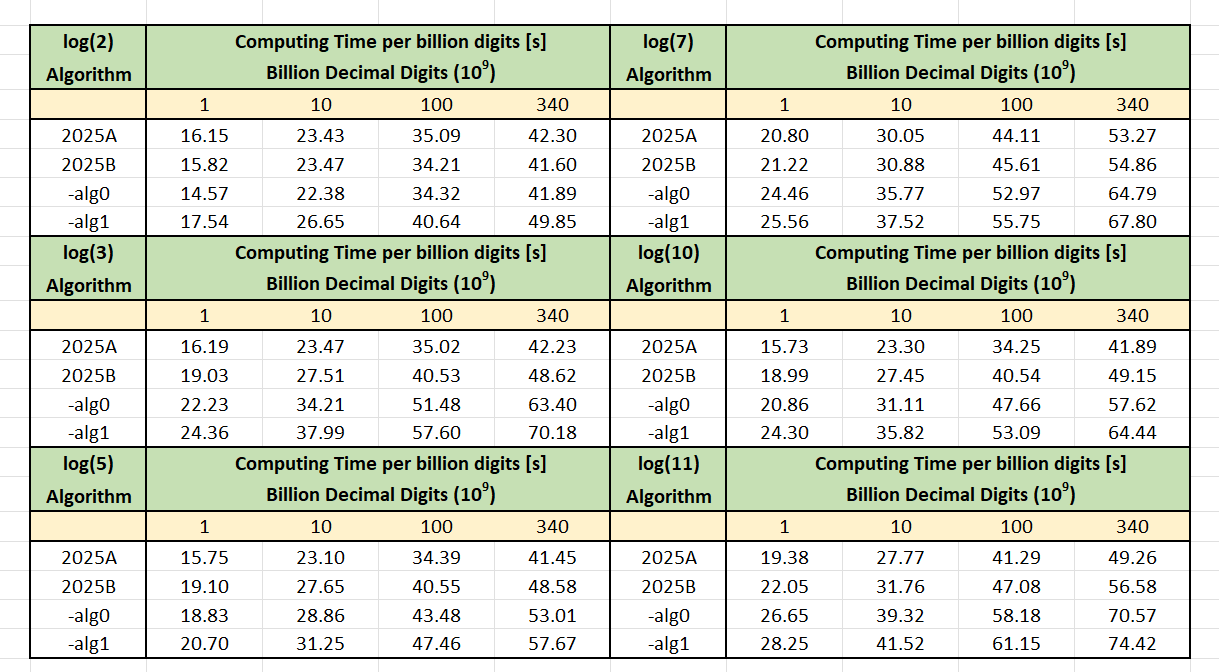}%
	\vspace{-6pt}
	\caption{Computing Time per billion digits in \texttt{[s]} for \texttt{log} algorithms.}
\end{figure}\vspace{0pt}

It is interesting to observe in Figure 2 how overall performance deteriorates and efficiency decays as more digits are calculated. This fact is associated to the specific y-cruncher software implementation. The widening performance gap between the different algorithms is also evident except for \texttt{log(2)} where algorithms \texttt{-alg0} and \texttt{2025B} perform almost the same at the highest level of decimal places. Note that Figure 2 gives a reciprocal of speed measurement, therefore lower values are better. The next page shows Figure 2 graphic performance comparisons between these algorithms.\newline\pagebreak
	\includegraphics[width=\textwidth]%
	{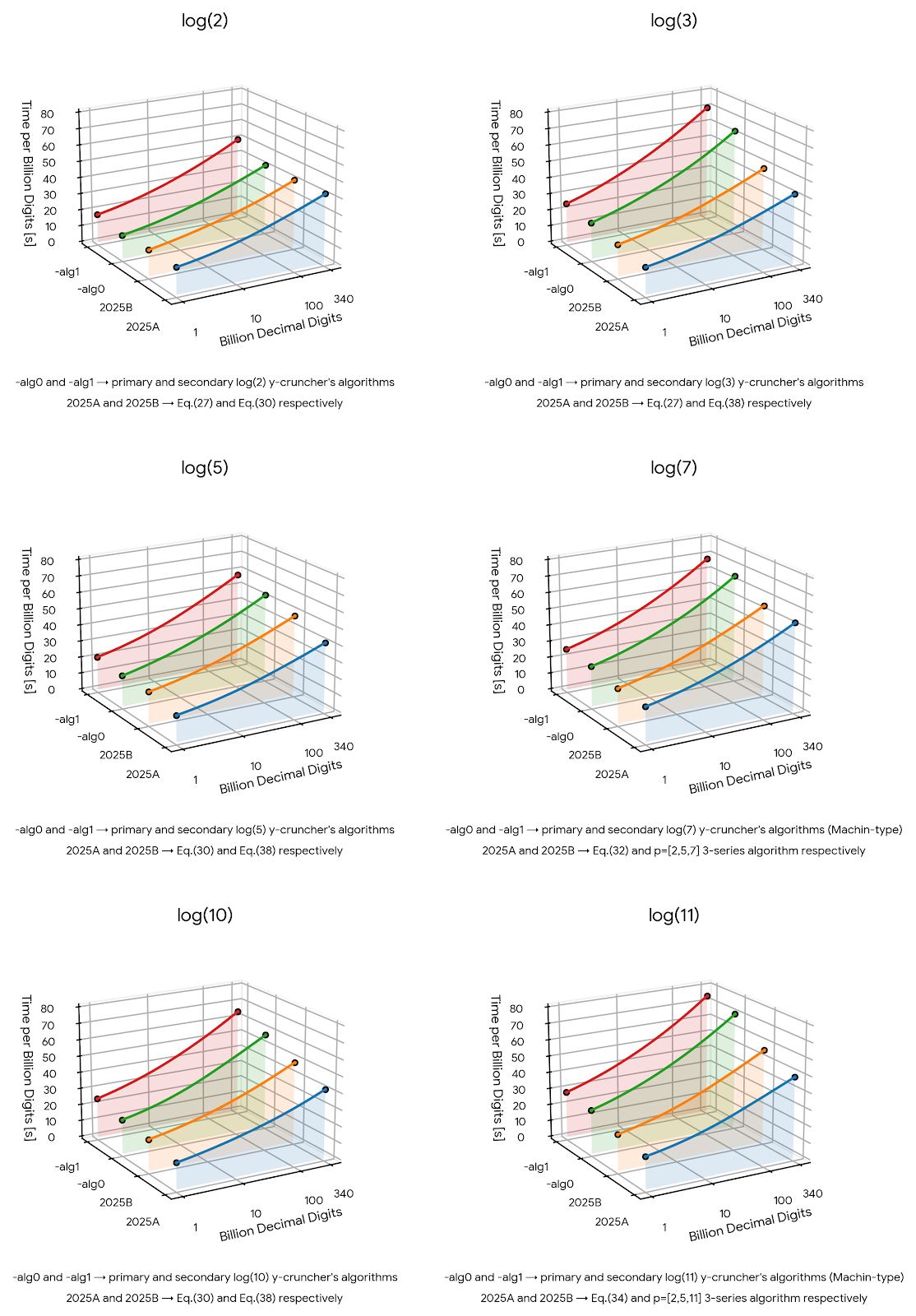}%
\pagebreak
\section{Fast series for logarithms of the sequence of primes}

\bigskip

A different application of this multi-valuation of logarithms algorithm is the
fast computing of some elementary functions. In particular the novel technique
introduced in \cite{FREDJ}, to get efficient multiprecision evaluations of
$\exp(x)$ and $\log(x)$ for $x$ real depend on the fast computing of the
sequence of logarithms of multiple primes. In this section several exponent
solutions $_{\text{\texttt{b}}}\boldsymbol{X}_{\boldsymbol{\pi}%
\text{\texttt{(n)}}}$ for\ different input \texttt{b-bit} sizes and \texttt{n =
5,\thinspace8,\thinspace11 }and \texttt{17} are presented.

\subsection{\textit{n} = 5}

$\boldsymbol{\pi}$\texttt{(5) = [2,3,5,7,11],} with input \texttt{b =
[64,96,128,192]-bit}, \ in these cases total binary splitting costs are
\texttt{C}$_{\text{s}}\ $\texttt{=\ [1.02940,0.778551,0.718053,0.572104]} and
costs per logarithm \texttt{[0.205880,0.155710,0.143611,0.114421]}
respectively, having these solution matrices with primes power exponents,
\[%
\begin{tabular}
[c]{clc}%
${}_{\text{\texttt{64}}}\boldsymbol{X}_{\boldsymbol{\pi}\text{\texttt{(5)}}} $
&  & $_{\text{\texttt{128}}}\boldsymbol{X}_{\boldsymbol{\pi}\text{\texttt{(5)}%
}}$\\
&  & \\
$\left[
\begin{tabular}
[c]{rrrrr}%
\texttt{-3} & \texttt{2} & \texttt{-1} & \texttt{2} & \texttt{-1}\\
\texttt{-7} & \texttt{-1} & \texttt{1} & \texttt{1} & \texttt{1}\\
\texttt{-1} & \texttt{5} & \texttt{0} & \texttt{0} & \texttt{-2}\\
\texttt{-4} & \texttt{0} & \texttt{2} & \texttt{1} & \texttt{-1}\\
\texttt{0} & \texttt{-5} & \texttt{1} & \texttt{2} & \texttt{0}%
\end{tabular}
\right]  $ &  & $\left[
\begin{tabular}
[c]{rrrrr}%
\texttt{-1} & \texttt{2} & \texttt{-4} & \texttt{5} & \texttt{-2}\\
\texttt{-3} & \texttt{4} & \texttt{-2} & \texttt{-2} & \texttt{2}\\
\texttt{-1} & \texttt{-7} & \texttt{4} & \texttt{1} & \texttt{0}\\
\texttt{-4} & \texttt{1} & \texttt{7} & \texttt{0} & \texttt{-4}\\
\texttt{-17} & \texttt{5} & \texttt{0} & \texttt{2} & \texttt{1}%
\end{tabular}
\right]  $\\
&  & \\
${}_{\text{\texttt{96}}}\boldsymbol{X}_{\boldsymbol{\pi}\text{\texttt{(5)}}} $
&  & $_{\text{\texttt{192}}}\boldsymbol{X}_{\boldsymbol{\pi}\text{\texttt{(5)}%
}}$\\
&  & \\
$\left[
\begin{tabular}
[c]{rrrrr}%
\texttt{-3} & \texttt{4} & \texttt{-2} & \texttt{-2} & \texttt{2}\\
\texttt{-1} & \texttt{-7} & \texttt{4} & \texttt{1} & \texttt{0}\\
\texttt{-5} & \texttt{-1} & \texttt{-2} & \texttt{4} & \texttt{0}\\
\texttt{-9} & \texttt{2} & \texttt{4} & \texttt{0} & \texttt{-1}\\
\texttt{-5} & \texttt{1} & \texttt{-3} & \texttt{0} & \texttt{3}%
\end{tabular}
\right]  $ &  & $\left[
\begin{tabular}
[c]{rrrrr}%
\texttt{-4} & \texttt{2} & \texttt{-11} & \texttt{2} & \texttt{6}\\
\texttt{-2} & \texttt{2} & \texttt{2} & \texttt{-7} & \texttt{4}\\
\texttt{-26} & \texttt{1} & \texttt{0} & \texttt{5} & \texttt{3}\\
\texttt{-1} & \texttt{13} & \texttt{3} & \texttt{-7} & \texttt{-2}\\
\texttt{-2} & \texttt{15} & \texttt{-1} & \texttt{-2} & \texttt{-4}%
\end{tabular}
\right]  $%
\end{tabular}
{}\medskip
\]
For \texttt{b = 64-bit}\ the mapped series column vector Eqs.(\ref{24}%
--\ref{25}) can be written as%
\[%
\genfrac{}{}{0pt}{0}{_{\text{\texttt{64}}}\boldsymbol{S}_{\boldsymbol{\pi
}\text{\texttt{(5)}}}\genfrac{}{}{0pt}{0}{\genfrac{}{}{0pt}{1}{{}}{{}}}{{}%
}}{\left[
\begin{tabular}
[c]{rrrrr}%
\texttt{4776624687651518} & \texttt{-796104114607999} & \texttt{2} &
\texttt{1} & \texttt{3156153406910380800}\\
\texttt{2420321086147582} & \texttt{-403386847690751} & \texttt{2} &
\texttt{1} & \texttt{1395919399595212800}\\
\texttt{241517233468190} & \texttt{-40252872244375} & \texttt{2} & \texttt{1}
& \texttt{87851769180634800}\\
\texttt{47948189888898} & \texttt{-7991364981249} & \texttt{-2} & \texttt{1} &
\texttt{12622326837120000}\\
\texttt{249074611495424} & \texttt{-41512435244032} & \texttt{64} & \texttt{1}
& \texttt{1424387459508300}%
\end{tabular}
\right]  }%
\vspace*{8pt}%
\]
and the linear combination coefficients%
\begin{gather*}
_{\text{\texttt{64}}}\boldsymbol{X}_{\boldsymbol{\pi}\text{\texttt{(5)}}}^{-1}%
\genfrac{}{}{0pt}{0}{\genfrac{}{}{0pt}{1}{{}}{{}}}{{}}%
\\
\left[
\begin{tabular}
[c]{rrrrr}%
\texttt{-24} & \texttt{27} & \texttt{46} & \texttt{-41} & \texttt{31}\\
\texttt{-38} & \texttt{43} & \texttt{73} & \texttt{-65} & \texttt{49}\\
\texttt{-56} & \texttt{63} & \texttt{107} & \texttt{-95} & \texttt{72}\\
\texttt{-67} & \texttt{76} & \texttt{129} & \texttt{-115} & \texttt{87}\\
\texttt{-83} & \texttt{94} & \texttt{159} & \texttt{-142} & \texttt{107}%
\end{tabular}
\right]
\end{gather*}

\subsection{\textit{n} = 8}

$\boldsymbol{\pi}$\texttt{(8) = [2,3,5,7,11,13,17,19],} with \texttt{b =
[64,96,128,192]-bit}, \ in these cases total binary splitting costs are
\texttt{C}$_{\text{s}}\ $\texttt{=\ [1.52370,1.02976,0.815916,0.714236]} and
costs per logarithm \texttt{[0.190463,0.128720,0.101989,0.0892794]}
respectively, giving%
\[%
\begin{tabular}
[c]{cc}%
${}_{\text{\texttt{64}}}\boldsymbol{X}_{\boldsymbol{\pi}\text{\texttt{(8)}}} $
& $_{\text{\texttt{128}}}\boldsymbol{X}_{\boldsymbol{\pi}\text{\texttt{(8)}}}%
$\\
& \\
\multicolumn{1}{l}{$\!\hspace*{-4pt}\left[
\begin{tabular}
[c]{rrrrrrrr}%
\texttt{-3} & \texttt{-1} & \texttt{1} & \texttt{1} & \texttt{0} & \texttt{1}
& \texttt{0} & \texttt{-1}\\
\texttt{-1} & \texttt{2} & \texttt{0} & \texttt{2} & \texttt{0} & \texttt{-1}
& \texttt{-1} & \texttt{0}\\
\texttt{-3} & \texttt{2} & \texttt{-1} & \texttt{2} & \texttt{-1} & \texttt{0}
& \texttt{0} & \texttt{0}\\
\texttt{-4} & \texttt{1} & \texttt{-2} & \texttt{1} & \texttt{0} & \texttt{0}
& \texttt{0} & \texttt{1}\\
\texttt{-7} & \texttt{-1} & \texttt{1} & \texttt{1} & \texttt{1} & \texttt{0}
& \texttt{0} & \texttt{0}\\
\texttt{-1} & \texttt{1} & \texttt{3} & \texttt{0} & \texttt{-1} & \texttt{0}
& \texttt{-1} & \texttt{0}\\
\texttt{-2} & \texttt{1} & \texttt{0} & \texttt{-1} & \texttt{2} & \texttt{-1}
& \texttt{0} & \texttt{0}\\
\texttt{-1} & \texttt{3} & \texttt{-2} & \texttt{-1} & \texttt{0} & \texttt{1}
& \texttt{0} & \texttt{0}%
\end{tabular}
\right]  $} & \multicolumn{1}{l}{$\hspace*{-4pt}\!\left[
\begin{tabular}
[c]{rrrrrrrr}%
\texttt{-10} & \texttt{-2} & \texttt{0} & \texttt{-1} & \texttt{-1} &
\texttt{3} & \texttt{1} & \texttt{1}\\
\texttt{-2} & \texttt{3} & \texttt{1} & \texttt{-1} & \texttt{-3} & \texttt{1}
& \texttt{-1} & \texttt{2}\\
\texttt{-1} & \texttt{7} & \texttt{2} & \texttt{-2} & \texttt{1} & \texttt{0}
& \texttt{-1} & \texttt{-2}\\
\texttt{-2} & \texttt{2} & \texttt{-1} & \texttt{-5} & \texttt{0} & \texttt{3}
& \texttt{1} & \texttt{0}\\
\texttt{-4} & \texttt{4} & \texttt{-1} & \texttt{4} & \texttt{-1} &
\texttt{-1} & \texttt{-1} & \texttt{0}\\
\texttt{-1} & \texttt{0} & \texttt{0} & \texttt{0} & \texttt{2} & \texttt{-5}
& \texttt{1} & \texttt{2}\\
\texttt{-5} & \texttt{-1} & \texttt{-1} & \texttt{-1} & \texttt{2} &
\texttt{-1} & \texttt{0} & \texttt{2}\\
\texttt{0} & \texttt{-1} & \texttt{-3} & \texttt{1} & \texttt{3} & \texttt{1}
& \texttt{-1} & \texttt{-1}%
\end{tabular}
\right]  $}\\
\multicolumn{1}{l}{} & \multicolumn{1}{l}{}\\
${}_{\text{\texttt{96}}}\boldsymbol{X}_{\boldsymbol{\pi}\text{\texttt{(8)}}} $
& $_{\text{\texttt{192}}}\boldsymbol{X}_{\boldsymbol{\pi}\text{\texttt{(8)}}}%
$\\
\multicolumn{1}{l}{} & \multicolumn{1}{l}{}\\
$\hspace*{-4pt}\!\left[
\begin{tabular}
[c]{rrrrrrrr}%
\texttt{-6} & \texttt{-2} & \texttt{-2} & \texttt{1} & \texttt{2} & \texttt{0}
& \texttt{1} & \texttt{0}\\
\texttt{-2} & \texttt{-3} & \texttt{1} & \texttt{-1} & \texttt{0} & \texttt{2}
& \texttt{1} & \texttt{-1}\\
\texttt{-8} & \texttt{1} & \texttt{-1} & \texttt{0} & \texttt{-1} & \texttt{1}
& \texttt{0} & \texttt{2}\\
\texttt{-6} & \texttt{1} & \texttt{0} & \texttt{3} & \texttt{-1} & \texttt{1}
& \texttt{0} & \texttt{-1}\\
\texttt{-2} & \texttt{-5} & \texttt{1} & \texttt{1} & \texttt{0} & \texttt{-1}
& \texttt{0} & \texttt{2}\\
\texttt{-3} & \texttt{2} & \texttt{3} & \texttt{-1} & \texttt{1} & \texttt{-1}
& \texttt{-1} & \texttt{0}\\
\texttt{-1} & \texttt{-1} & \texttt{-1} & \texttt{2} & \texttt{0} & \texttt{1}
& \texttt{1} & \texttt{-2}\\
\texttt{-3} & \texttt{0} & \texttt{3} & \texttt{1} & \texttt{-1} & \texttt{-2}
& \texttt{1} & \texttt{0}%
\end{tabular}
\right]  $ & $\hspace*{-4pt}\!\left[
\begin{tabular}
[c]{rrrrrrrr}%
\texttt{-1} & \texttt{-4} & \texttt{-1} & \texttt{1} & \texttt{-4} &
\texttt{1} & \texttt{0} & \texttt{4}\\
\texttt{-8} & \texttt{-5} & \texttt{-1} & \texttt{0} & \texttt{2} & \texttt{2}
& \texttt{2} & \texttt{-1}\\
\texttt{-2} & \texttt{2} & \texttt{-7} & \texttt{1} & \texttt{1} & \texttt{0}
& \texttt{-2} & \texttt{4}\\
\texttt{-10} & \texttt{-2} & \texttt{0} & \texttt{-1} & \texttt{-1} &
\texttt{3} & \texttt{1} & \texttt{1}\\
\texttt{-4} & \texttt{-3} & \texttt{3} & \texttt{0} & \texttt{-1} &
\texttt{-3} & \texttt{4} & \texttt{0}\\
\texttt{-6} & \texttt{-5} & \texttt{6} & \texttt{-3} & \texttt{1} &
\texttt{-1} & \texttt{-1} & \texttt{3}\\
\texttt{-15} & \texttt{9} & \texttt{1} & \texttt{1} & \texttt{-1} & \texttt{2}
& \texttt{-1} & \texttt{-1}\\
\texttt{-4} & \texttt{-1} & \texttt{-2} & \texttt{4} & \texttt{0} & \texttt{3}
& \texttt{-4} & \texttt{1}%
\end{tabular}
\right]  $%
\end{tabular}
\vspace*{4pt}{}%
\]
For \texttt{b = 64-bit}, mapped series $_{\text{\texttt{64}}}\boldsymbol{S}%
_{\boldsymbol{\pi}\text{\texttt{(8)}}}$--Eq.(\ref{25})-- and linear
combination coefficients $_{\text{\texttt{64}}}\boldsymbol{X}_{\boldsymbol{\pi
}\text{\texttt{(8)}}}^{-1}$\ are\vspace*{-4pt}%
\[%
\genfrac{}{}{0pt}{0}{_{\text{\texttt{64}}}\boldsymbol{S}_{\boldsymbol{\pi
}\text{\texttt{(8)}}}\genfrac{}{}{0pt}{0}{\genfrac{}{}{0pt}{1}{{}}{{}}}{{}%
}}{\left[
\begin{tabular}
[c]{rrrrr}%
\texttt{5647204593443138} & \texttt{-941200765573249} & \texttt{-2} &
\texttt{1} & \texttt{3858450988743187200}\\
\texttt{4831089675144514} & \texttt{-805181612523497} & \texttt{-2} &
\texttt{1} & \texttt{3199387769558770608}\\
\texttt{4776624687651518} & \texttt{-796104114607999} & \texttt{2} &
\texttt{1} & \texttt{3156153406910380800}\\
\texttt{2930726371203202} & \texttt{-488454395200001} & \texttt{-2} &
\texttt{1} & \texttt{1756236860945280000}\\
\texttt{2420321086147582} & \texttt{-403386847690751} & \texttt{2} &
\texttt{1} & \texttt{1395919399595212800}\\
\texttt{2121535521565502} & \texttt{-353589253593751} & \texttt{2} &
\texttt{1} & \texttt{1191771871116750000}\\
\texttt{1827734856225826} & \texttt{-304622476037153} & \texttt{-2} &
\texttt{1} & \texttt{996571801837579968}\\
\texttt{1523460247052798} & \texttt{-253910041174999} & \texttt{2} &
\texttt{1} & \texttt{800958681571230000}%
\end{tabular}
\right]  }%
\]
\vspace*{-12pt}
\begin{gather*}
_{\text{\texttt{64}}}\boldsymbol{X}_{\boldsymbol{\pi}\text{\texttt{(8)}}}^{-1}%
\genfrac{}{}{0pt}{0}{\genfrac{}{}{0pt}{1}{{}}{{}}}{{}}%
\\
\left[
\begin{tabular}
[c]{rrrrrrrr}%
\texttt{-72} & \texttt{26} & \texttt{-4} & \texttt{-72} & \texttt{84} &
\texttt{-26} & \texttt{-57} & \texttt{41}\\
\texttt{-114} & \texttt{41} & \texttt{-6} & \texttt{-114} & \texttt{133} &
\texttt{-41} & \texttt{-90} & \texttt{65}\\
\texttt{-167} & \texttt{60} & \texttt{-9} & \texttt{-167} & \texttt{195} &
\texttt{-60} & \texttt{-132} & \texttt{95}\\
\texttt{-202} & \texttt{73} & \texttt{-11} & \texttt{-202} & \texttt{236} &
\texttt{-73} & \texttt{-160} & \texttt{115}\\
\texttt{-249} & \texttt{90} & \texttt{-14} & \texttt{-249} & \texttt{291} &
\texttt{-90} & \texttt{-197} & \texttt{142}\\
\texttt{-266} & \texttt{96} & \texttt{-15} & \texttt{-266} & \texttt{311} &
\texttt{-96} & \texttt{-211} & \texttt{152}\\
\texttt{-294} & \texttt{105} & \texttt{-15} & \texttt{-294} & \texttt{343} &
\texttt{-106} & \texttt{-232} & \texttt{167}\\
\texttt{-306} & \texttt{110} & \texttt{-17} & \texttt{-305} & \texttt{357} &
\texttt{-110} & \texttt{-242} & \texttt{174}%
\end{tabular}
\right]
\end{gather*}

\subsection{\textit{n} = 11}

$\boldsymbol{\pi}$\texttt{(11) = [2,3,5,7,11,13,17,19,23,29,31]} and \texttt{b
= [128,256]-bit},\ in these cases total binary splitting costs are
\texttt{C}$_{\text{s}}\ $\texttt{=\ [1.04932,0.674650]} and costs per
logarithm \texttt{[0.0953929,0.0613318]} respectively. The next matrices
having the following exponents of primes, were obtained.%
\[
\]
Using \texttt{SearchLLLog(primes(11),128,-9,256,20000,1743454858545)}, this
solution is found for input \texttt{b = 128-bit}%

\begin{gather*}
_{\text{\texttt{128}}}\boldsymbol{X}_{\boldsymbol{\pi}\text{\texttt{(11)}}}%
\genfrac{}{}{0pt}{0}{\genfrac{}{}{0pt}{1}{{}}{{}}}{{}}%
\\
\left[
\begin{tabular}
[c]{rrrrrrrrrrr}%
\texttt{-10} & \texttt{-2} & \texttt{0} & \texttt{-1} & \texttt{-1} &
\texttt{3} & \texttt{1} & \texttt{1} & \texttt{0} & \texttt{0} & \texttt{0}\\
\texttt{-2} & \texttt{3} & \texttt{1} & \texttt{-1} & \texttt{-3} & \texttt{1}
& \texttt{-1} & \texttt{2} & \texttt{0} & \texttt{0} & \texttt{0}\\
\texttt{-5} & \texttt{6} & \texttt{0} & \texttt{-2} & \texttt{0} & \texttt{0}
& \texttt{-1} & \texttt{0} & \texttt{-1} & \texttt{2} & \texttt{0}\\
\texttt{-1} & \texttt{7} & \texttt{2} & \texttt{-2} & \texttt{1} & \texttt{0}
& \texttt{-1} & \texttt{-2} & \texttt{0} & \texttt{0} & \texttt{0}\\
\texttt{0} & \texttt{0} & \texttt{-5} & \texttt{3} & \texttt{1} & \texttt{1} &
\texttt{0} & \texttt{-2} & \texttt{1} & \texttt{0} & \texttt{0}\\
\texttt{-2} & \texttt{-2} & \texttt{1} & \texttt{-2} & \texttt{1} & \texttt{1}
& \texttt{-2} & \texttt{0} & \texttt{1} & \texttt{0} & \texttt{1}\\
\texttt{-2} & \texttt{0} & \texttt{0} & \texttt{-1} & \texttt{1} & \texttt{3}
& \texttt{0} & \texttt{1} & \texttt{-2} & \texttt{0} & \texttt{-1}\\
\texttt{-8} & \texttt{2} & \texttt{3} & \texttt{-1} & \texttt{-1} & \texttt{1}
& \texttt{0} & \texttt{0} & \texttt{-1} & \texttt{0} & \texttt{1}\\
\texttt{\ -3} & \texttt{2} & \texttt{-1} & \texttt{1} & \texttt{0} &
\texttt{0} & \texttt{2} & \texttt{-2} & \texttt{1} & \texttt{-1} &
\texttt{0}\\
\texttt{-2} & \texttt{2} & \texttt{-1} & \texttt{-5} & \texttt{0} & \texttt{3}
& \texttt{1} & \texttt{0} & \texttt{0} & \texttt{0} & \texttt{0}\\
\texttt{-1} & \texttt{0} & \texttt{0} & \texttt{0} & \texttt{2} & \texttt{-5}
& \texttt{1} & \texttt{2} & \texttt{0} & \texttt{0} & \texttt{0}%
\end{tabular}
\right]
\end{gather*}
\vspace*{2pt}%
\[
\]
By means of \texttt{SearchLLLog(primes(11),256,-13,384,20000,1743448649013)},
this solution is obtained for input \texttt{b = 256-bit}
\begin{gather*}
_{\text{\texttt{256}}}\boldsymbol{X}_{\boldsymbol{\pi}\text{\texttt{(11)}}}%
\genfrac{}{}{0pt}{0}{\genfrac{}{}{0pt}{1}{{}}{{}}}{{}}%
\\
\left[
\begin{tabular}
[c]{rrrrrrrrrrr}%
\texttt{2} & \texttt{5} & \texttt{-2} & \texttt{5} & \texttt{-1} & \texttt{0}
& \texttt{1} & \texttt{1} & \texttt{0} & \texttt{-6} & \texttt{1}\\
\texttt{14} & \texttt{-4} & \texttt{-7} & \texttt{-2} & \texttt{1} &
\texttt{0} & \texttt{-3} & \texttt{1} & \texttt{2} & \texttt{2} & \texttt{0}\\
\texttt{13} & \texttt{-1} & \texttt{3} & \texttt{-6} & \texttt{1} &
\texttt{-2} & \texttt{0} & \texttt{3} & \texttt{1} & \texttt{0} &
\texttt{-3}\\
\texttt{2} & \texttt{6} & \texttt{2} & \texttt{-4} & \texttt{-1} & \texttt{-2}
& \texttt{0} & \texttt{-2} & \texttt{1} & \texttt{0} & \texttt{2}\\
\texttt{12} & \texttt{10} & \texttt{-5} & \texttt{2} & \texttt{-2} &
\texttt{2} & \texttt{-1} & \texttt{-1} & \texttt{-2} & \texttt{0} &
\texttt{-1}\\
\texttt{3} & \texttt{-9} & \texttt{1} & \texttt{-1} & \texttt{-2} & \texttt{0}
& \texttt{2} & \texttt{1} & \texttt{-3} & \texttt{0} & \texttt{4}\\
\texttt{\ 11} & \texttt{-7} & \texttt{0} & \texttt{5} & \texttt{1} &
\texttt{0} & \texttt{-2} & \texttt{1} & \texttt{-3} & \texttt{-1} &
\texttt{1}\\
\texttt{6} & \texttt{-8} & \texttt{-6} & \texttt{0} & \texttt{0} & \texttt{0}
& \texttt{6} & \texttt{-1} & \texttt{-1} & \texttt{1} & \texttt{0}\\
\texttt{1} & \texttt{5} & \texttt{-3} & \texttt{-6} & \texttt{6} & \texttt{-3}
& \texttt{0} & \texttt{-1} & \texttt{1} & \texttt{0} & \texttt{1}\\
\texttt{1} & \texttt{0} & \texttt{0} & \texttt{8} & \texttt{-4} & \texttt{1} &
\texttt{-2} & \texttt{0} & \texttt{0} & \texttt{2} & \texttt{-3}\\
\texttt{1} & \texttt{-1} & \texttt{0} & \texttt{0} & \texttt{-1} & \texttt{9}
& \texttt{-2} & \texttt{-4} & \texttt{-2} & \texttt{0} & \texttt{1}%
\end{tabular}
\right]
\end{gather*}%
\[
\]
Note that last parameter in argument list of function \texttt{SearchLLLog} is
the seed used to get these solution matrices. This is set to repeat
experiments. However \texttt{m = 20000} (\texttt{nmax} parameter) is large
enough to capture at random the same solution with different seeds for both
\texttt{b = [128,256]-bit} cases, so this \texttt{seed} can be left blank on
\texttt{SearchLLLog} calling. As the parameter \texttt{n} increases, the
probability of identifying the optimal solution decreases. However, the
sub-optimal solutions discovered by the LLL--Monte Carlo algorithm still yield
highly efficient sets of series, which are sufficient to achieve the required
numerical objectives.

\subsection{\textit{n} = 17}

$\boldsymbol{\pi}$\texttt{(17)\thinspace=\thinspace\lbrack
2,3,5,7,11,13,17,19,23,29,31,37,41,43,47,53,59]}, with input \texttt{b = 256-bit}.
By using \texttt{SearchLLLog(primes(17),256,-14,384,20000)} it gives a binary
splitting cost \texttt{C}$_{\text{s}}\ $\texttt{=\ 0.846356} and cost per
logarithm\ \texttt{0.0497857} with initial random \texttt{seed =
1743461930745}. This provides the following solution matrix%
\begin{gather*}
_{\text{\texttt{256}}}\boldsymbol{X}_{\boldsymbol{\pi}\text{\texttt{(17)}}}%
\genfrac{}{}{0pt}{0}{\genfrac{}{}{0pt}{1}{{}}{{}}}{{}}%
\\
\left[
\begin{tabular}
[c]{rrrrrrrrrrrrrrrrr}%
\texttt{3} & \texttt{1} & \texttt{-2} & \texttt{-6} & \texttt{2} & \texttt{2}
& \texttt{0} & \texttt{-2} & \texttt{0} & \texttt{1} & \texttt{0} &
\texttt{-2} & \texttt{1} & \texttt{0} & \texttt{1} & \texttt{1} & \texttt{0}\\
\texttt{1} & \texttt{-5} & \texttt{2} & \texttt{0} & \texttt{-1} & \texttt{-2}
& \texttt{1} & \texttt{-1} & \texttt{2} & \texttt{0} & \texttt{1} & \texttt{1}
& \texttt{0} & \texttt{-1} & \texttt{1} & \texttt{1} & \texttt{-2}\\
\texttt{1} & \texttt{1} & \texttt{-4} & \texttt{-1} & \texttt{2} & \texttt{0}
& \texttt{-1} & \texttt{-2} & \texttt{1} & \texttt{2} & \texttt{-2} &
\texttt{0} & \texttt{2} & \texttt{-1} & \texttt{1} & \texttt{0} & \texttt{0}\\
\texttt{6} & \texttt{6} & \texttt{2} & \texttt{0} & \texttt{1} & \texttt{-2} &
\texttt{0} & \texttt{0} & \texttt{0} & \texttt{-2} & \texttt{0} & \texttt{1} &
\texttt{1} & \texttt{-2} & \texttt{1} & \texttt{0} & \texttt{-2}\\
\texttt{0} & \texttt{0} & \texttt{7} & \texttt{2} & \texttt{0} & \texttt{-1} &
\texttt{2} & \texttt{-1} & \texttt{0} & \texttt{-1} & \texttt{-1} & \texttt{0}
& \texttt{1} & \texttt{1} & \texttt{0} & \texttt{-3} & \texttt{-1}\\
\texttt{3} & \texttt{5} & \texttt{1} & \texttt{0} & \texttt{0} & \texttt{4} &
\texttt{-1} & \texttt{-1} & \texttt{-2} & \texttt{0} & \texttt{-1} &
\texttt{1} & \texttt{1} & \texttt{-3} & \texttt{0} & \texttt{0} & \texttt{0}\\
\texttt{1} & \texttt{0} & \texttt{2} & \texttt{-2} & \texttt{0} & \texttt{0} &
\texttt{1} & \texttt{2} & \texttt{-5} & \texttt{1} & \texttt{0} & \texttt{2} &
\texttt{0} & \texttt{-1} & \texttt{0} & \texttt{-1} & \texttt{1}\\
\texttt{3} & \texttt{-1} & \texttt{6} & \texttt{0} & \texttt{-2} & \texttt{-1}
& \texttt{1} & \texttt{0} & \texttt{0} & \texttt{0} & \texttt{-3} &
\texttt{-2} & \texttt{1} & \texttt{0} & \texttt{2} & \texttt{0} & \texttt{0}\\
\texttt{0} & \texttt{-1} & \texttt{0} & \texttt{-2} & \texttt{-1} &
\texttt{-1} & \texttt{0} & \texttt{4} & \texttt{-1} & \texttt{-3} &
\texttt{-1} & \texttt{1} & \texttt{1} & \texttt{2} & \texttt{0} & \texttt{0} &
\texttt{0}\\
\texttt{1} & \texttt{-5} & \texttt{1} & \texttt{1} & \texttt{0} & \texttt{2} &
\texttt{0} & \texttt{0} & \texttt{-1} & \texttt{0} & \texttt{1} & \texttt{0} &
\texttt{3} & \texttt{-1} & \texttt{-2} & \texttt{-2} & \texttt{1}\\
\texttt{1} & \texttt{2} & \texttt{0} & \texttt{-1} & \texttt{-7} & \texttt{3}
& \texttt{-1} & \texttt{1} & \texttt{0} & \texttt{-1} & \texttt{0} &
\texttt{1} & \texttt{1} & \texttt{0} & \texttt{0} & \texttt{0} & \texttt{1}\\
\texttt{4} & \texttt{2} & \texttt{-2} & \texttt{-3} & \texttt{-1} &
\texttt{-1} & \texttt{2} & \texttt{0} & \texttt{2} & \texttt{-1} & \texttt{0}
& \texttt{-1} & \texttt{0} & \texttt{0} & \texttt{-1} & \texttt{2} &
\texttt{0}\\
\texttt{0} & \texttt{0} & \texttt{5} & \texttt{-2} & \texttt{0} & \texttt{4} &
\texttt{0} & \texttt{-2} & \texttt{0} & \texttt{0} & \texttt{0} & \texttt{2} &
\texttt{0} & \texttt{0} & \texttt{-2} & \texttt{-1} & \texttt{-1}\\
\texttt{0} & \texttt{-7} & \texttt{0} & \texttt{2} & \texttt{1} & \texttt{0} &
\texttt{0} & \texttt{-2} & \texttt{0} & \texttt{-1} & \texttt{3} & \texttt{-2}
& \texttt{1} & \texttt{0} & \texttt{0} & \texttt{2} & \texttt{-1}\\
\texttt{5} & \texttt{-9} & \texttt{-1} & \texttt{0} & \texttt{-3} & \texttt{1}
& \texttt{0} & \texttt{-2} & \texttt{1} & \texttt{0} & \texttt{-1} &
\texttt{0} & \texttt{1} & \texttt{3} & \texttt{1} & \texttt{0} & \texttt{0}\\
\texttt{3} & \texttt{-1} & \texttt{1} & \texttt{0} & \texttt{1} & \texttt{4} &
\texttt{1} & \texttt{-1} & \texttt{-1} & \texttt{1} & \texttt{0} & \texttt{1}
& \texttt{-3} & \texttt{-1} & \texttt{0} & \texttt{0} & \texttt{-1}\\
\texttt{2} & \texttt{0} & \texttt{0} & \texttt{-2} & \texttt{-3} & \texttt{2}
& \texttt{1} & \texttt{2} & \texttt{0} & \texttt{0} & \texttt{0} & \texttt{-1}
& \texttt{1} & \texttt{2} & \texttt{-2} & \texttt{0} & \texttt{-1}%
\end{tabular}
\right]
\end{gather*}\vspace*{-5pt}

\section{PARI GP Modules\bigskip}

Since there is no option of accompanying ancillary files beyond the source TeX file, 
PARI GP scripts \texttt{SearchLog} and \texttt{SearchLLLog} have been embedded as 
a large comment towards the end of the parent TeX file. This file shall be downloaded 
to extract (copy-paste) the code. Directions to install and test the scripts are placed at the 
beginning of such comment.
\smallskip
\section{Conclusions\bigskip}

This work has yielded new methodologies in the efficient computation of
logarithms and related functions through the development of a merged application of Ramanujan-type hypergeometric series, multi-valuation techniques and optimization.

\subsubsection{Ramanujan-Type Hypergeometric Series for Logarithms and
Arctangents}

Highly convergent Ramanujan-type hypergeometric series for $\log(x)$ as
$x\rightarrow1$, were derived leading to new series for a$\tanh(x)$ and
a$\tan(x)$. The rapid convergence of these series for small $|\,x\,|$ makes
them ideal for multiprecision computations. Furthermore, formulas being implemented in binary
splitting form were provided, which have potential applications in:\vspace*{3pt}

\quad$\circ$ Machin formulas for $\pi$

\quad$\circ$ Machin--like formulas for logarithms.

\quad$\circ$ Efficient computation of logarithms of integer sequences.\medskip

\subsubsection{Multi-Valuation Algorithms for Logarithm Computation}

New algorithms, based on very fast series for $\log(x)$ where $x\in%
\mathbb{Q}
$ and $x\rightarrow1$, were developed to determine optimal multi-valuation of
$n$ formulas for computing several logarithms simultaneously. Two effective
methods$^{1}$\footnotetext[1]{Only sequential feasible MILP solutions searching was performed (random seeks, like MORSE method \cite{MORSE}, were not implemented). The direct application of MILP solvers were not as efficient as the brute force method and therefore they were excluded from this conclusions section.} were carried out:\vspace*{3pt}

\quad$\circ$ A fast brute-force approach for $n<10$.

\quad$\circ$ A randomized Monte Carlo lattice reduction-based LLL method for
$n<50$.\vspace*{4pt}

These algorithms, coded as PARI/GP scripts, are embedded into the companion TeX file -- see Section 7.--, providing
efficient solutions to get optimal multi-valuation formulas.\vspace
*{-2pt}

\subsubsection{Applications}

The developed methods have been successfully applied to:\medskip

\textbf{\quad}$\circ$\textbf{\ }Fastest Series for Single Logarithms:\vspace
*{3pt}

\quad$\cdot$ The fastest known multiprecision series for log(2), log(3), log(5), log(7)
log(10) and log(11) constants were obtained using bi and tri-valuation techniques.

\quad$\cdot$ Other highly efficient and fastest series for single logarithms
of small natural numbers were also found, significantly expanding upon
the single series results of Part\daurl{:I}{https://arxiv.org/pdf/2506.08245}\aurl{:I}{}.

\quad$\cdot$ A specialized method for fast multiprecision computations of
$\log(n)$ for natural numbers $n$ was introduced, utilizing 4-valuation
formulas based on integers \texttt{[2,3,5,n]}. This method requires a small,
easily constructed storage table. Preliminary results show superior
performance compared to Machin-like formulas based on Luca-Najman
tables.\medskip

\textbf{\quad}$\circ$\textbf{\ }Highly Efficient Multi-Valuation Formulas for
Primes:\vspace*{3pt}

\quad$\cdot$ Novel tables for the fastest series of multi--valuated logarithms
of prime sequences were generated. These formulas are particularly valuable
for improving the efficiency of multi-prime argument reduction in elementary
functions computations.

%
%
%
\section{Acknowledgements\bigskip}

My deep acknowledgments to Jesús Guillera$^\dagger\,$(6 Feb 2026) for his inspiration and 
encouragement leading me to look for and find these unique hypergeometric formulas. To Henri Cohen, Bill Allombert, Karim Belabas and all the PARI GP team, without their 
highly efficient and wonderful computing tool there would have been very difficult to get these identities. Special acknowledgements to Alex Yee, for his lightning fast y-cruncher platform. 
Special thanks to Lorenz Milla, for his excellent disposition to test some of these 
new formulas extending the number of digits known to \texttt{2e12} decimal places  
for \texttt{log(10)} on \texttt{June 06 2025} with his powerful multicore setup by using bi-series Eq.(\ref{30}) for computing and 
tri-series Eq.(\ref{38}) for digits validation. (An efficiency improvement over \texttt{5x} in 
\texttt{4.75} years), I am in debt as well for his kind support to validate the edition of this paper. 
\\
\\
My deep thanks and full credits to Dmitriy Grigoryev as well, who kindly provided a dazzling fast 2TB RAM Intel Xeon 60--core setup generating the source data beyond $10^{11}$ decimal places to create the comparison tests tables and graphs located in Section 5, the main improvements in this 2026 version of the paper, without his support, it could not be done.
Special thanks also to Joshua Swanson and Bruce Sagan for motivating me to make public these works after walking 
for a while through MSE -- see \cite{JZ}-- and MO blogs.
\\

\pagebreak

\pagebreak 
\section{Appendix\bigskip}

The set of pairs $\{(p_j,x_j)_{j=1}^n\}$ in Eq.(\ref{12}) with $x_j\in\mathbb{Z}$ can be splitted as
\begin{equation}
	\{(p_j,x_j)_{j=1}^n\}=\{(p_j,x_j^+)_{j=1}^d\}\cup\,\{(p_j,-x_j^-)_{j=d+1}^n\}\tag{A.1}\label{A.1}
\end{equation}
 for some $d\in\mathbb{N}$ and $d \le n$, where $x^{+}=x$ if $x\ge0$ and $x^{-}=-x$ if $x<0$, so that for $p>0$ and $p\in\mathbb{Q}_{\ne1}$
\begin{equation}
p=\frac{u}{v}=\prod_{j=1}^{n}p_{j}^{x_{j}}\tag{A.2}\label{A.2}
\end{equation}
has numerator/denominator 
\begin{equation}
u=\prod_{j=1}^{d}p_{j}^{x_{j}^{+}}\,\ \text{and}\ \ v=\prod_{j=d+1}^{n}p_{j}^{x_{j}^{-}}\tag{A.3}\label{A.3}%
\end{equation}
respectively. Once $u,v$ have been found, these identities for parameters $\alpha,\beta,\gamma,\nu,\delta\,\in\,\mathbb{Z}_{\ne0}$ in Eq.(\ref{24}) defining the Ramanujan-type hypergeometric logarithm series,

\begin{equation}
	\mathcal{S}=\frac{1}{\gamma}\cdot%
	{\displaystyle\sum\limits_{k=1}^{\infty}}
	\frac{\alpha\hspace{1pt}k+\beta}{k(2k-1)}\cdot\left(  \frac
	{\mathcal{\nu}}{\mathcal{\delta}}\right)  ^{k}\cdot%
	\begin{bmatrix}%
		\genfrac{}{}{0pt}{0}{{}}{{}}%
		1 & \frac{1}{2}%
		\genfrac{}{}{0pt}{0}{{}}{{}}%
		\\%
		\genfrac{}{}{0pt}{0}{{}}{{}}%
		\frac{1}{6} & \frac{5}{6}%
		\genfrac{}{}{0pt}{0}{{}}{{}}%
	\end{bmatrix}
	_{k}\tag{A.4}\label{A.4}%
\end{equation}
are obtained from Eq.(\ref{9}),

\begin{equation}
	\alpha=-2\,(u+v)(u^{2}-14uv+v^{2})(u^{2}+4uv+v^{2})
	\tag{A.5}\label{A.5}%
\end{equation}

\begin{equation}
	\beta=(u+v)^{3}(u^{2}-8uv+v^{2})
	\tag{A.6}\label{A.6}%
\end{equation}

\begin{equation}
	\gamma=2\,(u-v)^{5}
	\tag{A.7}\label{A.7}%
\end{equation}

\begin{equation}
	\nu=(u-v)^{6}
	\tag{A.8}\label{A.8}%
\end{equation}

\begin{equation}
	\delta=108\ u^{2}\,v^{2}\,(u+v)^{2}
	\tag{A.9}\label{A.9}\vspace*{8pt}%
\end{equation}

Finally, simplified $(\alpha,\beta,\gamma)$ and $(\nu,\delta)$ are found by eliminating common factors

\begin{equation}
	(\alpha,\beta,\gamma)\ \longleftarrow\ \frac{(\alpha,\beta,\gamma)}{\gcd(\alpha,\beta,\gamma)}\tag{A.10}\label{A.10}
\end{equation}

\begin{equation}
	(\nu,\delta)\ \longleftarrow\ \frac{(\nu,\delta)}{\gcd(\nu,\delta)}\tag{A.11}\label{A.11}\vspace*{8pt}
\end{equation}

These parameters are applied on Eq.(\ref{A.4}) and Eqs.(\ref{24}--\ref{25}) to build the column vector of series $\boldsymbol{\,_bS_q}$.

\vfill\eject

\end{document}